\documentclass[10pt]{article}
\usepackage{amsthm}
\usepackage{amsfonts}
\usepackage{amsmath}
\usepackage{amssymb}
\usepackage[all]{xy}

\everymath{\displaystyle}
\CompileMatrices

\newtheorem{thm}{Theorem}[section]
\newtheorem{lemma}[thm]{Lemma}
\newtheorem{cor}[thm]{Corollary}
\newtheorem{prop}[thm]{Proposition}
\newtheorem{conj}[thm]{Conjecture}

\theoremstyle{definition}
\newtheorem{defn}[thm]{Definition}

\def\ve{\varepsilon}

\def\R{{\mathbb R}}
\def\Q{{\mathbb Q}}
\def\C{{\mathbb C}}

\def\Z{{\mathbb Z}}

\def\so{_{\scriptscriptstyle O}}
\def\su{_{\scriptscriptstyle U}}
\def\st{_{\scriptscriptstyle T}}

\def\po{^{\scriptscriptstyle O}}
\def\pu{^{\scriptscriptstyle U}}
\def\pt{^{\scriptscriptstyle T}}
\def\sr{_{\scriptscriptstyle \R}}
\def\sc{_{\scriptscriptstyle \C}}
\def\crt{^{\scriptscriptstyle {\it CRT}}}
\def\scrt{_{\scriptscriptstyle {\it CRT}}}
\def\ct{{\it CRT}}

\def\coker{\text {coker}}

\def\ext{\text {Ext}}
\def\hom{\text {Hom}}

\def\conj2#1{\setbox0=\hbox{$#1$}
\vbox{\hsize=\wd0\noindent\offinterlineskip
\hbox to\wd0{\cleaders \hbox{$\symb$}\hss}
\vskip1pt
\box0}}
\def\symb{\scriptscriptstyle\mathrel\sim\mkern-5mu}

\newcommand{\smv}[2]{       \bigl( \begin{smallmatrix} 
					{#1} \\ {#2}
                   		    \end{smallmatrix} \bigr)         }

\newcommand{\smh}[2]{       \bigl( \begin{smallmatrix} 
					{#1} & {#2}
                   		    \end{smallmatrix} \bigr)         }

\begin{document}
\bibliographystyle{amsplain}


\title{Real C*-Algebras, United $KK$-Theory, \\ and the Universal Coefficient
Theorem} \date{\today}
\author{Jeffrey L. Boersema}

\maketitle

\vspace*{-0.525in}
\hspace*{4.71in}

\begin{center}
\quad \\[-10pt]
Department of Mathematics\\ 
Seattle University\\
Seattle, WA  98133\\
{\sf boersema@seattleu.edu} 
\end{center}

\setcounter{page}{1}

\thispagestyle{empty}

\vspace{0.2in}

\begin{abstract}  
\baselineskip15pt

We define united $KK$-theory for real C*-algebras $A$ and $B$ such that $A$
is separable and $B$ is $\sigma$-unital, extending united $K$-theory in the
sense that $KK\crt(\R, B) = K\crt(B)$.  United $KK$-theory contains real,
complex, and self-conjugate $KK$-theory; but unlike unaugmented real
$KK$-theory, it admits a universal coefficient theorem.  For all separable
$A$ and $B$ in which the complexification of $A$ is in the bootstrap
category, $KK\crt(A,B)$ can be written as the middle term of a short exact
sequence whose outer terms involve the united $K$-theory of $A$ and $B$. 
As a corollary, we prove that united $K$-theory classifies $KK$-equivalence
for real C*-algebras whose complexification is in the bootstrap category.
\end{abstract}

\bigskip
\bigskip

\baselineskip18pt


\section{Introduction}
The Universal Coefficient Theorem for the $K$-theory of complex
C*-algebras, proven by Rosenberg and Schochet (\cite{RS}) in 1987, states
the existence of an unnaturally split, short exact sequence 
$$0 \rightarrow
{\rm Ext}_{K_*(\C)}(K_*(A), K_*(B)) \xrightarrow{\kappa} KK\sc(A,B)
\xrightarrow{\gamma} {\rm Hom}_{K_*(\C)}(K_*(A), K_*(B)) \rightarrow 0$$ 
for any complex separable C*-algebras $A$ and $B$ such that $A$ is in the bootstrap
category $\mathcal{N}$.  Recall that $\mathcal{N}$ is the smallest
subcategory of complex, separable, nuclear C*-algebras which contains the
separable type I C*-algebras; which is closed under the operations of
taking inductive limits, stable isomorphisms, and crossed products by $\Z$
and $\R$; and which satisfies the two out of three rule for short exact
sequences (i.e. if $0 \rightarrow A \rightarrow B \rightarrow C \rightarrow
0$ is exact and two of $A$, $B$, $C$ are in $\mathcal{N}$, then the third
is also in $\mathcal{N})$.
  
One immediate and powerful corollary of this theorem states that
$KK$-equivalence is completely characterized by $K$-theory for algebras in
the bootstrap class, $\mathcal{N}$ (see Section 7 of~\cite{RS},
Section~23.10 of \cite{Black}, and Section~2.4 of \cite{Rordam} for this
and other corollaries).  The role of the Universal Coefficient Theorem in
the standard toolkit in the subject of $K$-theory of complex C*-algebras is
considerable.  In particular, it is used essentially in the proofs of
$K$-theoretic classification theorems for simple C*-algebras, such as those
found in \cite{Kirchberg} and \cite{Phillips} in the purely infinite case; in
\cite{EG} in the stably finite, real rank zero case; and in
\cite{Lin} in the tracial rank zero case.  

In this paper, we consider real C*-algebras.  Recall that a real 
C*-algebra $A$ is a Banach $*$-algebra over $\R$ which satisfies the 
C*-equation $\|a^* a\| = \|a\|^2$ as well as the axiom that $1 + a^*a$ is 
invertible in the unitization $\widetilde{A}$ for all $a \in A$.  This is 
equivalent (by \cite{Palmer}; see also Chapter 15 of \cite{Goodearl}) to 
saying that $A$ is isomorphic to a 
norm-closed adjoint-closed algebra of operators on a Hilbert space over 
$\R$.  Note that every complex C*-algebra can be considered a real C*-algebra by 
forgetting the complex structure.  Conversely, given any real C*-algebra 
$A$, we can form the complexification $\C \otimes A = A\sc$.

Except for the classification of complex AF-algebras in \cite{Elliot},
which was repeated for real AF-algebras in \cite{Gio}, none of the major
work in classifying complex C*-algebras has been carried over to real
C*-algebras.  One reason for this is that there is no universal coefficient
theorem for real C*-algebras.  While Kasparov \cite{Kasp} considered
$KK$-theory in a very general setting, Rosenberg and Schochet proved the
universal coefficient theorem only for complex C*-algebras as they
explained, ``For reasons pointed out already by Atiyah, there can be no
good K\"unneth Theorem or Universal Coefficient Theorem for the KKO groups
of real C*-algebras; this explains why we deal only with complex
C*-algebras.''  The reference to Atiyah is \cite{Ati2} wherein a
counter-example is given; the essential obstruction lies in the
homological algebra, namely that $K_*(A)$ does not have projective or
injective dimension one in the category of $K_*(\R)$-modules.  In fact, 
the projective dimension may be infinite.

In the present paper, we remedy this situation by proving a universal
coefficient theorem for real C*-algebras using united $K$-theory and 
united $KK$-theory.  This paper is a sequel to
our earlier paper \cite{Boer2} in which we extended united $K$-theory
from its original topological setting \cite{Bou} to the setting of real
C*-algebras.  In the present paper, we further extend these constructions,
defining united $KK$-theory for real C*-algebras.  That is, we define
$KK\crt(A,B)$ for real C*-algebras $A$ and $B$ such that $A$ is separable
and $B$ is $\sigma$-unital.  This functor combines real $KK$-theory,
complex $KK$-theory, and self-conjugate $KK$-theory as well as the natural
transformations among the three.  Thus we will show that it takes values in
the category of so-called \ct-modules (also referred to as the category
\ct), which is a much more hospitable setting than the category of
$K_*(\R)$-modules (the main source of information on the category \ct~and
its hospitality is \cite{Bou}).  In \cite{Boer2}, we took advantage of the
fact that the object $K\crt(A)$ has projective dimension one in \ct~to
prove a K\"unneth formula for real C*-algebras.  In the current paper, we
take advantage of the fact that the same object has injective dimension one
to prove the following universal coefficient theorem.
\begin{thm}[Main theorem]
Let $A$ and $B$ be real separable C*-algebras 
such that $\C \otimes A$ is in the bootstrap category $\mathcal{N}$.  
Then there is a short
exact sequence
\end{thm}
$$0 \rightarrow {\rm Ext}\scrt(K\crt(A), K\crt(B)) \xrightarrow{\kappa}
	KK\crt(A,B) \xrightarrow{\gamma}
	{\rm Hom}\scrt(K\crt(A), K\crt(B)) 
	\rightarrow 0$$
where $\kappa$ has degree $-1$.
As in the K\"unneth sequence for united $K$-theory, this short exact 
sequence does not split in general.

If Rosenberg and Schochet say that no good UCT is possible for real
C*-algebras, we believe that ours is the best possible.  Indeed our
universal coefficient theorem does exactly what a universal coefficient
theorem is supposed to do: it expresses $KK\crt(A,B)$ (and thus $KK_*(A,
B)$, since united $KK$-theory contains real $KK$-theory) in terms of
$K$-theoretic data involving $A$ and $B$; although the data required is
more than just $K_*(A)$ and $K_*(B)$.

Furthermore, as a corollary to our universal coefficient theorem, we find
that two real separable C*-algebras whose complexifications are in the bootstrap
category are $KK$-equivalent if and only if their united $K$-theories are
isomorphic.  If the development of new mathematical machinery is justified
by its ability to answer questions that exist independently of that
machinery, then this corollary justifies the development of united
$K$-theory.  Indeed, united $K$-theory exactly captures $KK$-equivalence, a
notion which involves only elements in real $KK$-theory and is thus prior
to the development of united $K$-theory or united $KK$-theory.  This also 
strengthens our case that united $K$-theory is the right functor to look at
for real C*-algebras, especially when considering the possibility of a
classification of real simple C*-algebras.  

The work of Hewitt in
\cite{Hewitt} on \ct-modules suggests that the same information 
can be retained in a somewhat smaller invariant.  On the other hand,
Bousfield has made us aware of examples of acyclic \ct-modules which are
not determined by their real part.  Together with our result that every
acyclic \ct-module can be realized as the united $K$-theory of a real
C*-algebra (which will be written and published elsewhere), this implies
that real $K$-theory alone is certainly not sufficient to determine
$KK$-equivalence.

Our proof of the universal coefficient theorem will be found in the final
section of the paper, and follows the same line of attack as in \cite{RS}. 
We first prove that $\gamma$ is an isomorphism in case $K\crt(B)$ is an
injective \ct-module.  We then tackle the general case by using a geometric
injective resolution of $K\crt(B)$.  The fact that $K\crt(B)$ has injective
dimension one is a necessary condition for the existence of such a
geometric injective resolution.

In the intervening sections we build up the required machinery.  In
Section~2, we define united $KK$-theory and work out its main properties. 
Section~3 is entirely algebraic, defining $\hom\scrt(M,N)$ in the category
of \ct-modules and working out its main properties.


\section{United $KK$-theory}
In this section, we define united $KK$-theory and the intersection product 
for united $KK$-theory.  In this paper, each $K$-group or $KK$-group is 
considered as $\Z$-graded module
with a periodicity automorphism of degree 
8 (or degree 2 if the C*-algebras are complex).  This periodicity 
automorphism is implicit in the $K_*(\R)$-module structure described below.  

\subsection{The United $KK$-Theory Bifunctor}

Recall from \cite{Boer2} that $T = \{f \colon [0,1] 
\rightarrow \C \mid f(0) = \overline{f(1)} \}$ is the self-conjugate 
algebra used in the definition of self-conjugate $K$-theory for real 
C*-algebras.  We will use it 
here to define self-conjugate $KK$-theory.  

\begin{defn} For any pair of real C*-algebras $A$ and $B$ such that $A$ is separable 
and $B$ is $\sigma$-unital, we define real $KK$-theory, complex 
$KK$-theory, and self-conjugate $KK$-theory as follows:
\begin{enumerate}
\item[(1)] $KKO_*(A,B) = KK_*(A,B)$
\item[(2)] $KKU_*(A,B) = KK_*(A, \C \otimes B)$
\item[(3)] $KKT_*(A,B) = KK_*(A, T \otimes B)$
\end{enumerate}
\end{defn}

Because of the intersection product, real $KK$-theory takes values in the 
category of modules over the ring 
$K_*(\R)$, given in degrees $0$ through $8$ by  
$$
K_*(\R) \hspace{.5cm} = \hspace{.5cm}
\Z \hspace{1cm} 
\Z_2 \hspace{1cm}
\Z_2 \hspace{1cm}
0 \hspace{1cm}
\Z \hspace{1cm}
0 \hspace{1cm}
0 \hspace{1cm}
0 \hspace{1cm}
\Z  
$$
(see page 23 in \cite{Schroder} or Table~1 in \cite{Boer2}).  The 
generators are the elements $\eta\so \in K_0(\R)$, $\xi \in K_4(\R)$, 
and the invertible element $\beta\so \in K_8(\R)$.
Complex $KK$-theory takes values in the category of modules over the ring 
$K_*(\C)$, which is the free polynomial ring generated by $\beta\su \in K_2(\C)$:
$$
K_*(\C) \hspace{.5cm} = \hspace{.5cm}
\Z \hspace{1cm} 
0 \hspace{1cm}
\Z \hspace{1cm}
0 \hspace{1cm}
\Z \hspace{1cm}
0 \hspace{1cm}
\Z \hspace{1cm}
0 \hspace{1cm}
\Z  $$
(see Table 2 in \cite{Boer2}).  The action of $\beta\su$ on $KK_*(A, \C 
\otimes B)$ is actually implemented by first letting $\beta\su$ pass through the 
natural homomorphisms
$$K_2(\C) = KK_2(\R, \C) \rightarrow KK_2(\C, \C \otimes \C) \rightarrow KK_2(\C, \C) \; $$
and then applying the intersection product.  We let $\beta\su$ also denote 
the element of $KK_2(\C, \C)$ which more directly implements the natural transformation 
on $KK_*(A, \C \otimes B)$.  
Finally, self conjugate $KK$-theory takes values in the category of modules 
over the ring $K_*(T) = KK_*(\R, 
T)$
$$
K_*(T) \hspace{.5cm} = \hspace{.5cm}
\Z \hspace{1cm} 
\Z_2 \hspace{1cm}
0 \hspace{1cm}
\Z \hspace{1cm}
\Z \hspace{1cm}
\Z_2 \hspace{1cm}
0 \hspace{1cm}
\Z \hspace{1cm}
\Z  $$
(see Table 3 in \cite{Boer2}) with generators
$\eta\st$ in degree 1, $\omega$ in degree 3, and the invertible element $\beta\st$ 
is degree 4.  As in the complex case, the action of these elements on 
$KK(A, T \otimes B)$ is implemented by intersection product with 
corresponding elements in $KK_*(T, T)$.

Furthermore, there are natural transformation 
\begin{align*}
c_n &\colon KKO_n(A, B) \rightarrow KKU_n(A, B) \\
r_n &\colon KKU_n(A, B) \rightarrow KKO_n(A, B) \\ 
\ve_n &\colon KKO_n(A, B) \rightarrow KKT_n(A, B) \\ 
\zeta_n &\colon KKT_n(A, B) \rightarrow KKU_n(A, B) \\
(\psi\su)_n &\colon KKU_n(A, B) \rightarrow KKU_n(A, B) \\
(\psi\st)_n &\colon KKT_n(A, B) \rightarrow KKT_n(A, B)\\ 
\gamma_n &\colon KKU_n(A, B) \rightarrow KKT_{n-1}(A, B) \\
\tau_n &\colon KKT_n(A, B) \rightarrow KKO_{n+1}(A, B) 
\end{align*}
among the three variants of 
$KK$-theory given by intersection product with the $KK$-elements
\begin{align*}
c &\in KK_0(\R, \C)&
r &\in KK_0(\C, \R)\\
\ve &\in KK_0(\R, T)&
\zeta &\in KK_0(T, \C)\\
\psi\su &\in KK_0(\C, \C)&
\psi\st &\in KK_0(T, T)\\
\gamma &\in KK_{-1}(\C, T)&
\tau &\in KK_1(T, \R)  \; .
\end{align*}

These elements are obtained as follows.
In Section~1 of \cite{Boer2}, the operations $c$, $r$, $\ve$, $\zeta$, $\psi\su$,
$\psi\st$, and $\gamma$ are defined on united $K$-theory.  Each is induced by a
homomorphism among the C*-algebras $\R$, $\C$, and $T$ (or suspensions thereof or
matrix algebras thereover).  Hence each produces a $KK$-element which 
implements the natural transformation.  Also in Section~1 of
\cite{Boer2}, the operation $\tau$ is defined in terms of two homomorphisms
$\sigma_i \colon T \rightarrow C(S^1, M_4(\R))$ (for $i=1, 2$) such that $\pi
\circ \sigma_1 = \pi \circ \sigma_2$ where $\pi \colon C(S^1, M_4(\R)) \rightarrow
M_4(\R)$ is the base-point evaluation map.  Therefore there are $KK$-elements
$\sigma_i \in KK_0(T, C(S^1, \R))$ such that $\sigma_1 - \sigma_2$ lies in the
kernel of the homomorphism $KK_0(T, C(S^1, \R)) \rightarrow KK_0(T, \R)$.  Hence
$\sigma_1 - \sigma_2$ lifts to produce a unique element $\tau \in KK_0(T, S\R)
= KK_1(T, \R)$.

\begin{defn} [United $KK$-Theory]
For any pair of real C*-algebras $A$ and $B$ such that $A$ is separable 
and $B$ is $\sigma$-unital, we define the united
$KK$-theory to the be the triple
$$KK\crt(A,B) := \{ KKO_*(A,B), KKU_*(A,B), KKT_*(A,B) \}\, $$
consisting of graded modules over $K_*(\R)$, $K_*(\C)$, and $K_*(T)$ 
(respectively)
together with the operations $\{c, r, \ve, \zeta, \psi\su, \psi\st, \gamma, 
\tau \}$.
\end{defn}

\begin{thm} \label{unitedproperties}
United $KK$-theory is a natural bifunctor; co-variant in the second argument
and contra-variant in the first argument; additive, homotopy invariant, 
stable, and has period 8 in both arguments; and it satisfies the following 
exactness properties:

Assume that $0 \longrightarrow 
A \xrightarrow{\alpha}  
B \xrightarrow{\beta}
C \longrightarrow 0$ is a semisplit short exact sequence of 
$\sigma$-unital C*-algebras
(semisplit means that $\beta$ has a completely positive cross-section with norm
no more than 1). If $D$ is separable, there is a long exact sequence
$$ \cdots \longrightarrow 
KK\crt(D,A) \xrightarrow{\alpha_*} 
KK\crt(D,B) \xrightarrow{\beta_*}
KK\crt(D,C) \xrightarrow{\delta} 
KK\crt(D,A) \longrightarrow \cdots \; .$$
If $B$ is separable, there 
is a long exact sequence 
$$ \cdots \longrightarrow 
KK\crt(A,D) \xleftarrow{\alpha^*} 
KK\crt(B,D) \xleftarrow{\beta^*}
KK\crt(C,D) \xleftarrow{\delta} 
KK\crt(A,D) \longleftarrow \cdots \; .$$
In either case, $\delta$ is a homomorphism of degree $-1$ implemented by 
an element $\delta \in KK_{-1}(C, A)$.
\end{thm}

\begin{proof}
All properties follow from the corresponding properties enjoyed by 
$KK$-theory for real C*-algebras (see Sections~2.3-2.5 in \cite{Schochet}). 

For exactness, assume that
\begin{equation} \label{Rsequence}
0 \rightarrow A \rightarrow B \rightarrow C \rightarrow 0
\end{equation}
is semisplit.  Recall Lemma~1.23 in \cite{Boer2} which states that a short exact
sequence is semisplit if and only if its complexification  
\begin{equation} \label{Csequence}
0 \rightarrow \C \otimes A \rightarrow \C \otimes B \rightarrow \C
\otimes C \rightarrow 0 \; 
\end{equation}
is also semisplit.  Let $s$ be a completely positive section of 
Extension~\ref{Csequence} with norm at
most 1.  If we apply the functor $C(S^1, -)$ to this last sequence we obtain
the sequence \begin{equation} \label{S^1sequence} 0 \rightarrow C(S^1,\C
\otimes A) \rightarrow C(S^1,\C \otimes B)  \rightarrow C(S^1,\C \otimes C)
\rightarrow 0 \;  \end{equation}
which is, we claim, semisplit with section $s_* \colon f \mapsto s \circ f$.  
It is easily seen that $\|s_*\| \leq 1$ since $\|s\| \leq 1$.  That $s_*$ is 
completely positive 
follows from the fact that a function $f \in C(X, D)$ is positive if and 
only if $f(x)$ is positive for all $x \in X$.  Therefore
Sequence~\ref{S^1sequence} is semisplit.  

Now recall that $T = \{f \in C(S^1, \C) \mid f(z) = \overline{f(-z)} \}$.  Then
$\C \otimes T = C(S^1, \C)$.  Hence Sequence~\ref{S^1sequence} is the
complexification of the sequence \begin{equation} \label{Tsequence}
0 \rightarrow T \otimes A \rightarrow T \otimes B \rightarrow T
\otimes C \rightarrow 0 \; .
\end{equation}
Therefore, by Lemma~1.23 in \cite{Boer2} again, Sequence~\ref{Tsequence} is
semisplit.

To summarize, if Sequence~\ref{Rsequence} is semisplit, so are
Sequences~\ref{Csequence} and \ref{Tsequence}.  Therefore by Theorem~2.5.6 in
\cite{Schroder}, there are long exact sequences in each of the three parts of
united $KK$-theory.  The maps in these long exact sequences commute with the
transformation maps of united $KK$-theory and therefore we have long exact
sequences in united $KK$-theory as in the statement of the Theorem.
\end{proof}

\begin{prop} \label{KKcrt}
If $A$ and $B$ are real C*-algebras such that $A$ is separable 
and $B$ is $\sigma$-unital, then $KK(A,B)$ is an acyclic \ct-module.
\end{prop}

The statement that $KK(A,B)$ is a \ct-module (see Section~2 of \cite{Bou} 
or Section~1.3 of \cite{Boer2}) means 
that $KK_*(A,B)$ is a module
over $KO_*(\R) = KK_*(\R, \R)$; that $KK_*(A, \C \otimes B)$ is a module over
$KU_*(\R) = KK_*(\R, \C)$; that $KK_*(A, T \otimes B)$ is a module over
$KT_*(\R) = KK_*(\R, T)$;
and that the following relations hold among the
operations:
\begin{align*}  
rc &= 2    & \psi\su \beta\su &= -\beta\su \psi\su & \xi &= r \beta\su^2 c \\
cr &= 1 + \psi\su & \psi\st \beta\st &= \beta\st \psi\st 
	&  \omega &= \beta\st \gamma \zeta \\ 
r &= \tau \gamma & \varepsilon \beta\so &= \beta\st^2 \varepsilon 
	& \beta\st \varepsilon \tau 
			&= \varepsilon \tau \beta\st + \eta\st \beta\st    \\
c &= \zeta \varepsilon & \zeta \beta\st &= \beta\su^2 \zeta
	 &  \varepsilon r \zeta &= 1 + \psi\st   \\
(\psi\su)^2 &= 1 & \gamma \beta\su^2 &= \beta\st \gamma       
	&  \gamma c \tau &= 1 - \psi\st \\
(\psi\st)^2 &= 1 & \tau \beta\st^2 &= \beta\so \tau & \tau &= -\tau \psi\st  \\
\psi\st \varepsilon &= \varepsilon & 
	\gamma &= \gamma \psi\su 
	\qquad & \tau \beta\st \varepsilon &= 0 \\
\zeta \gamma &= 0 & \eta\so &= \tau \varepsilon 
	& \varepsilon \xi &= 2 \beta\st \varepsilon \\
\zeta &= \psi\su \zeta & \eta\st &= \gamma \beta\su \zeta 
	& \xi \tau &= 2 \tau \beta\st \; .
\end{align*}

The statement that $KK(A,B)$ is acyclic means that the following sequences are
exact:
\begin{equation} \label{ouseq} 
\cdots \longrightarrow  KK_n(A,B) \xrightarrow{\eta\so} 
KK_{n+1}(A,B) \xrightarrow{c} 
KK_{n+1}(A, \C \otimes B) \xrightarrow{r \beta\su^{-1}}
KK_{n-1}(A,B) \longrightarrow \cdots \end{equation}
\begin{equation} \label{otseq} 
\cdots \longrightarrow  KK_{n}(A,B) \xrightarrow{\eta\so^2} 
KK_{n+2}(A,B) \xrightarrow{\varepsilon} 
KK_{n+2}(A, T \otimes B) \xrightarrow{\tau \beta\st^{-1}}
KK_{n-1}(A,B) \longrightarrow \cdots \end{equation}
\begin{equation} \label{utseq} 
\cdots \rightarrow  KK_{n+1}(A, \C \otimes B) \xrightarrow{\gamma} 
KK_n(A,T \otimes B) \xrightarrow{\zeta} 
KK_n(A, \C \otimes B) \xrightarrow{1 - \psi\su}
KK_{n}(A, \C \otimes B) \rightarrow \cdots \end{equation}

\begin{proof}[Proof of Proposition~\ref{KKcrt}]
We first show that $KK(A,B)$ is acyclic.  
In Section~1.4 of \cite{Boer2} we
produced two semisplit short exact sequences which were homotopy equivalent to
\begin{equation} \label{rc}
0 \rightarrow S^{-1}\R \rightarrow \R \xrightarrow{c} \C \rightarrow 0    
\end{equation}
and
\begin{equation} \label{rt}
0 \rightarrow S^{-2}\R \rightarrow \R \xrightarrow{\ve} T \rightarrow 0 \;.
\end{equation}

In the first extension, 
it was proven that the homomorphism $S^{-1}\R \rightarrow \R$ is
represented by the $KK$-element $\eta\so \in KK_1(\R, \R)$ and the connecting 
homomorphism in the associated long exact sequence in $K$-theory
is represented by $r \beta\su^{-1} \in KK_{-2}(\C, \R)$.  For
the second it was proven that the homomorphism $S^{-2}\R \rightarrow \R$ is
represented by $\eta\so^2 \in KK_2(\R, \R)$ and the connecting homomorphism
is represented by $\tau \beta\st^{-1} \in KK_{-3}(T, \R)$. 

When we tensor these sequences by $B$ and apply part (5) of
Theorem~\ref{unitedproperties} we obtain Sequences~\ref{ouseq} and \ref{otseq}
as desired.

Sequence~\ref{utseq} is similarly derived from the split exact sequence
$$0 \rightarrow S\C \xrightarrow{\gamma} T 
\xrightarrow{\zeta} \C \rightarrow 0$$ 
where $\gamma$ is defined by inclusion and $\zeta$ is defined by evaluation
at one endpoint.  (See Theorem~1.5 of \cite{Boer2}.) 

Now we show that $KK(A,B)$ is a \ct-module.
We must show that the CRT-relations hold in this context.  Let $M$ and $N$ 
be elements of $\{\R, \C, T\}$.
Each relation is of the 
form $\psi = \phi$ where $\psi$ and $\phi$ are natural transformations from
$KK_*(A, M \otimes B)$ to $KK_*(A, N \otimes B)$ formed by taking the product
of one or more of the \ct-operations.  We will do better than show that $\psi$
and $\phi$ are equivalent transformations --- we will show that the
corresponding classes in $KK_*(M,N)$ are identical.

By Proposition~1.9 in \cite{Boer2} the \ct-relations hold for
united $K$-theory.  To show that they hold on the level of $KK$-elements it is
enough to show that the homomorphism $$\alpha \colon KK_*(M,N)
\rightarrow \hom_{K_*(\R)}(K_*(M), K_*(N))$$ is injective when $M$ and $N$ are
in $\{\R, \C, T\}$.  

We first need to compute $KK_*(M,N)$ for all possible combinations of $M$ and
$N$.  If $M = \R$, we refer to Table~1
in \cite{Boer2} since $KK_*(\R,N) \cong K_*(N)$.  We determine $KK_*(\C, \R)$
and $KK_*(T,\R)$ using the sequences
$$
KK_{n+1}(\R,\R) \xleftarrow{\eta\so} 
KK_n(\R,\R) \xleftarrow{c} 
KK_n(\C,\R) \xleftarrow{r \beta\su^{-1}} 
KK_{n+2}(\R, \R) \xleftarrow{\eta\so^{2}} 
KK_{n+1}(\R,\R) $$
and
$$ KK_{n+2}(\R, \R) \xleftarrow{\eta\so^2} 
KK_{n}(\R, \R) \xleftarrow{\ve}
KK_{n}(T, \R) \xleftarrow{\tau \beta\st^{-1}}
KK_{n+3}(\R, \R) \xleftarrow{\eta\so^2}
KK_{n+1}(\R, \R)$$
which are derived from \ref{rc} and \ref{rt}.
The remaining $KK$ groups are derived by using long exact 
sequences (again based on \ref{rc} and \ref{rt}) relating $KK(M,N)$ to $KK(\R, N)$ 
for $M,N \in \{\C, T\}$.  The results of these computations are recorded in 
Table~\ref{somekkgroups}

\begin{table}  
\caption{Some $KK$-groups} \label{somekkgroups} 
$$\begin{array}{|c|c|c|c|c|c|c|c|c|c|}  
\hline  \hline 
n & \makebox[1cm][c]{0} & \makebox[1cm][c]{1} & 
\makebox[1cm][c]{2} & \makebox[1cm][c]{3} 
& \makebox[1cm][c]{4} & \makebox[1cm][c]{5} 
& \makebox[1cm][c]{6} & \makebox[1cm][c]{7} 
& \makebox[1cm][c]{8} \\
\hline  \hline 
KK_*(\R,\R)  &
\Z & \Z_2 & \Z_2 & 0 & 
\Z & 0 & 0 & 0 & 
\Z         \\
\hline  
KK_*(\R, \C)   &
\Z & 0 & 
\Z & 0 & 
\Z & 0 & 
\Z & 0 & 
\Z         \\
\hline  
KK_*(\R, T)   &
\Z & \Z_2 & 0 & \Z &  
\Z & \Z_2 & 0 & \Z &
\Z         \\
\hline
KK_*(\C,\R)  &
\Z & 0 & \Z & 0 & 
\Z & 0 & \Z & 0 & 
\Z         \\
\hline  
KK_*(\C, \C)   &
\Z \oplus \Z & 0 & 
\Z \oplus \Z & 0 & 
\Z \oplus \Z & 0 & 
\Z \oplus \Z & 0 & 
\Z \oplus \Z         \\
\hline  
KK_*(\C, T)   &
\Z & \Z & \Z & \Z &  
\Z & \Z & \Z & \Z &
\Z         \\
\hline
KK_*(T,\R)  &
\Z & \Z & \Z_2 & 0 & 
\Z & \Z & \Z_2 & 0 & 
\Z         \\
\hline  
KK_*(T, \C)   &
\Z & \Z  & 
\Z & \Z & 
\Z & \Z  & 
\Z & \Z & 
\Z         \\
\hline
KK_*(T, T)   &
\Z \oplus \Z & \Z \oplus \Z_2 & 
\Z_2 & \Z & 
\Z \oplus \Z & \Z \oplus \Z_2 & 
\Z_2 & \Z & 
\Z \oplus \Z        \\
\hline \hline
\end{array}$$
\end{table}

Now we're ready to show that $\alpha$ is injective.  If $M=\R$, then for any
$\theta \in KK_*(M,N)$, we have  $\alpha(\theta) (1\so) = 1\so \cdot \theta =
\theta$.  Therefore, $\alpha$ is an isomorphism.

Now consider $KK_*(\C, \R)$.  The class $r \in KK_0(\C, \R)$ induces a non-trivial 
homomorphism of infinite degree from $K_{0}(\C) \cong \Z$ to $K_0(\R) \cong \Z$. 
Therefore the image of $\alpha$ in $\hom_{K_*(\R)}(K_*(\C), K_*(\R))$ contains an
infinite cyclic subgroup in degree $0$.  Since $KK_0(\C, \R) \cong \Z$, it
follows that $\alpha$ must be injective in degree 0.  More generally, the class $r
\beta\su^n \in KK_{2n}(\C, \R)$ induces a non-trivial homomorphism from
$K_{0}(\C)$ to $K_{2n}(\R)$ showing that $\alpha$ must be injective on
$KK_{2n}(\C, \R)$.  

According to Table~1 in \cite{Boer2}, $\psi\su \colon K_2(\C) \rightarrow
K_2(\C)$ is multiplication by $-1$.  Therefore the elements $\beta\su^n$ and
$\psi\su \beta\su^n$ induce two independent homomorphisms from $K_*(\C)$ to
itself, which show that $\alpha$ is injective for $KK_{2n}(\C, \C)$.

The rest of the proof proceeds in the same way.  Using Table~1 in \cite{Boer2}
we study the behavior of the operations among the $K$-theory of $\R$, $\C$, and
$T$.  The elements  
$\varepsilon r \beta\su^n$ and $\gamma \beta\su^n$ show that $\alpha$ is
injective for $KK_*(\C, T)$; $r \zeta \beta\st^n $, $\tau \beta\st^n$ , and
$\eta\so \tau \beta\st^n$ for $KK_*(T, \R)$;  $\beta\su^n \zeta$ and 
$\beta\su^n c \tau$ for $KK_*(T, \C)$; and finally, the elements $\beta\st^n$, $\psi\st
\beta\st^n$, $\eta\st \beta\st^n$, $\varepsilon \tau \beta\st^n$, $\eta\st
\varepsilon \tau \beta\st^n$, and $\omega \beta\st^n$ for $KK_*(T,T)$.     
\end{proof}


\subsection{The Intersection Product} \label{intproduct}

In this section we will describe an intersection product for united 
$KK$-theory and we will prove that this product is a {\it CRT}-pairing.  Throughout this section, we 
assume that all C*-algebras are separable and that enough of the 
C*-algebras are nuclear so that tensor products are unique.

The intersection product (see Section~2.4 of \cite{Schroder}) for
real $KK$-theory is a pairing 
$$\alpha \colon 
KK_i(A_1, B_1 \otimes D) \otimes KK_j(D \otimes A_2, B_2) \longrightarrow 
		KK_{i+j}(A_1 \otimes A_2, B_1 \otimes B_2)$$
which is natural, bilinear, and associative.  The associativity implies that the 
pairing is bilinear over the ring $KO_*(\R) = KK_*(\R, \R)$.  Hence we have a
$KO_*(\R)$-module pairing 
$$\alpha\so \colon KKO_*(A_1, B_1 \otimes D)
\otimes_{KO_*(\R)} KKO_*(D \otimes A_2, B_2) \longrightarrow  		
KKO_*(A_1\otimes A_2, B_1 \otimes B_2) \; . $$

We also have a $KU_*(\R)$-module pairing
$$\alpha\su \colon KKU_*(A_1, B_1 \otimes D) \otimes_{KU_*(\R)} 
KKU_*(D \otimes A_2, B_2) \longrightarrow 
		KKU_*(A_1 \otimes A_2, B_1 \otimes B_2)  $$
which arises by composing the intersection product with multiplication $\C \otimes\sr \C \rightarrow \C$.

Similarly, since $T$ is commutative, the homomorphism $T \otimes\sr T \rightarrow 
T$ inducess a $KT_*(\R)$-module 
pairing
$$\alpha\st \colon KKT_*(A_1, B_1 \otimes D) \otimes_{KT_*(\R)} 
KKT_*(D \otimes A_2, B_2) \longrightarrow KKT_*(A_1 \otimes A_2, B_1 \otimes B_2) 
	\; . $$

\begin{prop}
The pairings $\alpha\so$, $\alpha\su$, and $\alpha\st$ together form a natural 
associative {\it CRT}-pairing
$$\alpha \colon KK\crt(A_1, B_1 \otimes D) \otimes\scrt KK\crt(D \otimes A_2, B_2) 
\rightarrow
KK\crt(A_1 \otimes A_2, B_1 \otimes B_2)$$
\end{prop}

\begin{proof}
We must show that the \ct-pairing relations of 
Definition~3.1 in \cite{Boer2} hold.  These relations are
\begin{enumerate}
\item[{\rm(1)}] 
	$c(m\so \cdot\so n\so) = c(m\so) \cdot\so c(n\so)$ 
\item[{\rm(2)}]
	$r(c(m\so) \cdot\su n\su) = m\so \cdot\so r(n\su)$ 
\item[{\rm(3)}] 
	$r(m\su \cdot\su c(n\so)) = r(m\su) \cdot\so n\so$ 
\item[{\rm(4)}] 
	$\varepsilon(m\so \cdot\so n\so) 
		= \varepsilon(m\so) \cdot\st \varepsilon(n\so)$ 
\item[{\rm(5)}]
	$\zeta(m\st \cdot\st n\st) = \zeta(m\st) \cdot\su \zeta(n\st)$ 
\item[{\rm(6)}]
	$\psi\su(m\su \cdot\su n\su) = \psi\su(m\su) \cdot\su \psi\su(n\su)$ 
\item[{\rm(7)}]
	$\psi\st(m\st \cdot\st n\st) = \psi\st(m\st) \cdot\st \psi\st(n\st)$ 
\item[{\rm(8)}]
	$\gamma(m\su \cdot\su \zeta(n\st)) = \gamma(m\su) \cdot\st n\st$ 
\item[{\rm(9)}]
	$\gamma(\zeta(m\st) \cdot\su n\su) 
		= (-1)^{|m\st|} m\st \cdot\st \gamma(n\su)$ 
\item[{\rm(10)}] 
	$\tau(m\st \cdot\st \varepsilon(n\so)) = \tau(m\st) \cdot\so n\so$ 
\item[{\rm(11)}]
	$\tau(\varepsilon(m\so) \cdot\st n\st) 
		= (-1)^{|m\so|} m\so \cdot\so \tau(n\st)$  
\item[{\rm(12)}] $\varepsilon \tau(m\st \cdot\st n\st) 
		= \varepsilon \tau(m\st) \cdot\st n\st 		
 		+ (-1)^{|m\st|} m\st \cdot\st \varepsilon \tau (n\st)   
		+ \eta\st(m\st \cdot\st n\st) \; .$ 
\end{enumerate}
where $\cdot\so$, $\cdot\su$, and $\cdot\st$ represent the product under the
pairings $\alpha\so$, $\alpha\su$, and $\alpha\st$ respectively.

The first eleven relations hold for intersection products in $KK$-theory for the same 
reason that they hold for external products in $K$-theory:  diagrams on the level
of C*-algebras which commute up to homotopy (see Lemmas~1.7 and 1.10 in
\cite{Boer2}).  For example, the diagram
$$ \xymatrix{
S\C \otimes T				
\ar[r]^{1 \otimes \zeta}		
\ar[d]_{\gamma \otimes 1}	
& S\C \otimes \C			
\ar[r]					
& S\C					
\ar[d]_\gamma		     \\
T \otimes T
\ar[rr]
&& T
			}$$
commutes up to homotopy, as we have seen in the proof of Lemma~1.7 of
\cite{Boer2}.  Putting  in suspensions, we have the diagram 
$$ \xymatrix{
S^m\C \otimes S^n T				
\ar[r]^{1 \otimes \zeta}		
\ar[d]_{\gamma \otimes 1}	
& S^m\C \otimes S^n\C			
\ar[r]					
& S^{m+n}\C					
\ar[d]_\gamma	\\							
S^{m-1}T \otimes S^n T
\ar[rr]
&& S^{m+n-1} T	}$$
which commutes up to homotopy (recall that our convention is that the homomorphism 
$\gamma \colon S^n \C \rightarrow S^{n-1} T$ uses the leftmost suspension) 
and establishes the relation
$\gamma(x \cdot\su \zeta(y)) = \gamma(x) \cdot\st y$ for 
$x \in KK_m(A_1, \C \otimes B_1 \otimes D)$ and $y \in KK_n(D \otimes A_2,
T \otimes B_1)$.

The relation $\gamma(\zeta(x) \cdot\su y) = (-1)^m x \cdot\st \gamma(y)$ for 
$x \in KKT_m(A_1, B_1 \otimes D)$ and $y \in KKU_n(D \otimes A_2, B_2)$ is 
established similarly starting from the diagram 
$$ \xymatrix{
T \otimes S\C
	\ar[r]^{\zeta \otimes 1}
	\ar[d]_{1 \otimes \gamma}
& \C \otimes S\C
	\ar[r]
& S\C
	\ar[d]_\gamma			\\
T \otimes T
	\ar[rr]
&&T					}$$
which commutes up to homotopy.  Putting in suspensions we get the diagram
$$ \xymatrix{
S^m T \otimes S^n\C
	\ar[r]^{\zeta \otimes 1}
	\ar[d]_{1 \otimes \gamma}
& S^m \C \otimes S^n \C
	\ar[r]
& S^{m+n} \C
	\ar[d]_\gamma			\\
S^m T \otimes S^{n-1} T
	\ar[rr]
&& S^{m+n-1} T					}$$
which commutes up to homotopy and up to a factor of $\iota^n$ where $\iota$ is the 
involution on $S^2$ which interchanges the two suspension factors.  This factor of
$\iota$ arises because $\gamma$ is required to always use the outermost
suspension.  Since $\iota$ induces multiplication by $-1$ on $KK$-theory, the
relation is established. 

The rest of the first eleven {\it CRT}-pairing relations can also be established in 
this way based directly on the commutative diagrams used in establishing
Lemmas~1.7 and 1.10 in \cite{Boer2}.  

For the twelfth relation 
\begin{equation} \label{relation12}
	\varepsilon \tau(x \cdot\st y) = 
		\varepsilon \tau(x) \cdot\st y 
		+ (-1)^n x \cdot\st \varepsilon \tau(y) 
		+ \eta\st(x \cdot\st y)
\end{equation} 
for $x \in KKT_m(A_1, B_1 \otimes D)$ and $y \in KKT_n(D \otimes
A_2, B_2)$ we take a different tack.  Let $\phi\st$ be the
$KK$-element corresponding to the C*-algebra homomorphism $T \otimes T
\rightarrow T$.
The operation $x \otimes y
\mapsto \ve \tau(x \cdot\st y)$ is performed by composing the product
homomorphism
$$KK_m(A_1, T \otimes B_1 \otimes D) \otimes KK_n(D \otimes A_1, T \otimes B_2)
	\rightarrow
		KK_{m+n}(A_1 \otimes A_2, T \otimes T \otimes B_1 \otimes B_2)$$
by the homomorphism
$$KK_{m+n}(A_1 \otimes A_2, T \otimes T \otimes B_1 \otimes B_2)
	\rightarrow
		KK_{m+n+1}(A_1 \otimes A_2, T \otimes B_1 \otimes B_2)$$
that is induced by the $KK$-element 
	$\ve \tau \phi\st \in KK_1(T \otimes T, T)$.

Meanwhile the right side of Equation~\ref{relation12} is the image of $x
\otimes y$ under an operation which is obtained in exactly the same way, except
that the relevant $KK$-element is 
$\phi\st(\ve \tau \otimes 1) + (-1)^n \phi\st(1 \otimes \ve \tau) 
	+ \eta\st \phi\st$.  We will show that these two $KK$-elements are
identical.

We know that these two operations agree in the case $A_1 = A_2 = B_1 = B_2 = D
= \R$ by Proposition~4.1 of \cite{Boer2}.  Therefore, it is enough to show that
the map 
$$\alpha \colon KK_*(T \otimes T, T) \rightarrow \hom_*(K_*(T \otimes
T),K_*(T))$$  is injective.  Lemma~\ref{TtimesT}, which follows, implies the
existence of a split exact sequence 
$$0 \rightarrow ST \rightarrow T \otimes T \rightarrow T \rightarrow 0$$
which allows us to make the decompositions  
\begin{align*}
KK_*(T \otimes T, T) &\cong KK_*(ST,T) \oplus KK_*(T,T) \\
\hom_*(K_*(T \otimes T), K_*(T)) &\cong 
		\hom_*(K_*(ST), K_*(T)) \oplus \hom_*(K_*(T), K_*(T)) \; .
\end{align*}
Since $\alpha$ respects this direct sum decomposition, it suffices to show that
the  map 
$$KK_*(T, T) \rightarrow \hom(K_*(T), K_*(T))$$
is injective.  But this fact was established in the proof of 
Proposition~\ref{KKcrt}.
\end{proof}

\begin{lemma} \label{TtimesT}
There is a C*-algebra isomorphism $T \otimes T \cong C(S^1, T)$.
\end{lemma}

\begin{proof}
Let $I = [0,1]$.  Then 
$T \cong \{f \colon I \rightarrow \C \mid f(0) = \overline{f(1)} \}$ 
so we have 
$$T \otimes T \cong \{f \colon I^2 \rightarrow \C \otimes \C \mid 		
	f(0,t) = \overline{f(1,t)} 	\text{~and~} 
	f(s,0) = \conj2{f(s,1)}  \}$$
where the straight bar represents conjugation in the first factor of $\C \otimes
\C$ and the wavy bar represents conjugation in the second factor.

Now, recall from \cite{Boer2} the isomorphism $\C \oplus \C \rightarrow \C
\otimes \C$ given by $(\lambda_1, \lambda_2) \mapsto \tfrac{1}{2}(\lambda_1 +
\lambda_2) \otimes 1 + \tfrac{1}{2}(\lambda_1 - \lambda_2) i \otimes i$.  Under
this isomorphism, conjugation in the first factor of $\C \otimes \C$
corresponds to the involution $(\lambda_1, \lambda_2) \mapsto
(\overline{\lambda_2}, \overline{\lambda_1})$ while conjugation in the second
factor corresponds to the involution $(\lambda_1, \lambda_2) \mapsto
(\lambda_2, \lambda_1)$.

Therefore, 
$$T \otimes T \cong \{f,g \colon I^2 \rightarrow \C \mid 
		f(0,t) = \overline{g(1,t)},
		g(0,t) = \overline{f(1,t)},
		f(s,0) = g(s,1), \text{~and~}
		g(s,0) = f(s,1) \} \; .$$

Next, we use the gluing operation $(f,g) \mapsto h$ defined by 
$$h(s,t) = \begin{cases}
			f(s,2t) \qquad &t \leq 1/2 \\
			g(s,2t-1) \qquad &t \geq 1/2
	\end{cases}$$
which gives us an isomorphism
$$T \otimes T \cong \{h \colon I^2 \rightarrow \C \mid
	h(s,0) = h(s,1) \text{~and~} 
		h(0,t) = \overline{h(1,t + \tfrac{1}{2})} \} \; $$
where the addition $t + \tfrac{1}{2}$ is done modulo $1$.

Finally, we untwist this algebra using the operation $h \mapsto h'$
defined by
$h'(s,t) = h(s,t+s/2)$ where, again, addition is modulo $1$.  Hence we
finally have
$$
T \otimes T \cong \{h \colon I^2 \rightarrow \C \mid
	h(s,0) = h(s,1) \text{~and~} 
		h(0,t) = \overline{h(1,t)} \}   
	 \cong C(S^1, T)  $$
\end{proof}


\section{The Hom functor in \ct}  

This section is entirely algebraic, wherein we define the functor $\hom\scrt(M,N)$ 
and work out the properties of this functor necessary for our 
development of the universal coefficient theorem.  Note that 
$\hom\scrt(M,N)$ does not refer to the graded group of 
\ct-morphisms from $M$ to $N$, which we denote $[M,N]_*$.  
Instead, we will (in Section~ \ref{homdef}) define $\hom\scrt(M,N)$ 
to be a \ct-module $\hom\scrt(M,N)$ such that the adjoint relationship
$$[L, \hom\scrt(M,N)]_* \cong [L \otimes\scrt M, N]_* \; $$
holds for all \ct-modules $L$, $M$, and $N$.  (See section 3 of
\cite{Boer2} for tensor products in \ct.)

In Section~\ref{homfree}, we prove three propositions 
which characterize $\hom\scrt(M,N)$ when $M$ is a
monogenic free $\ct$-module.  These theorems help us get a more concrete handle 
on the \hom functor in \ct.  Indeed, based on these theorems, one strategy to 
compute $\hom\scrt(M,N)$ when $M$ is acyclic is to build a free resolution of 
$M$ and apply the functor $\hom\scrt(-, N)$ to the resolution.  Finally, in 
Section~\ref{homacyc}, 
we prove a couple propositions giving sufficient conditions for acyclicity of 
$\hom\scrt(M,N)$.  

\subsection{Defining $\hom\scrt(M,N)$} \label{homdef}

In what follows, for any $i \in \Z$ and any $X \in \{O,U,T\}$, let 
$F(b,i,X)$ be the free \ct-module with a single generator $b \in 
F(b,i,X)^X_i$.  These modules are described in Section~2.4 of \cite{Bou} and 
in Section~2.1 of \cite{Boer2}.

\begin{defn}
Given two \ct-modules $M$ and $N$, we define the \ct-module $\hom\scrt(M,N)$ as
follows
\begin{enumerate}
\item[(1)] For $i \in \Z$ and $X \in \{O,U,T\}$, we define $\hom\scrt(M,N)_i^X$ to
be the group of \ct-module homomorphisms of degree 0 from $F(b,i,X) \otimes\scrt M$
to $N$.  That is,  	
	$$\hom\scrt(M,N)_i^X = [F(b,i,X) \otimes\scrt M, N]_0$$
\item[(2)] Let $\theta \colon L_i^X \rightarrow L_j^Y$ be a \ct-operation 
and let $\mu_\theta$ be the function $F(c,j,Y) \rightarrow F(b,i,X)$
which carries $c$ to $\theta(b) \in F(b,i,X)_j^Y$.  Then for any $\alpha 
\in \hom\scrt(M,N)_i^X = [F(b,i,X) \otimes M,N]_0$, define $\theta(\alpha) 
\in \hom\scrt(M,N)_j^Y = [F(c,j,Y) \otimes\scrt M,N]_0$ by the formula 
$\theta(\alpha) = \alpha \circ(\mu_\theta \otimes 1)$.  This describes the 
operation
$$\theta \colon \hom\scrt(M,N)_i^X \rightarrow \hom\scrt(M,N)_j^Y$$
\end{enumerate}
\end{defn}

\begin{prop}
$\hom\scrt(-,-)$ is a bifunctor with values in \ct~which is
contravariant in the first argument and covariant in the second.

\end{prop}

\begin{proof}
The only statement that is not clear is that $\hom\scrt(M,N)$ is a \ct-module if $M$
and $N$ are.  To show that the \ct-relations are all satisfied, it suffices to show
that if $\phi \circ \psi = \theta$ is a general \ct-relation, then it
holds in $\hom\scrt(M,N)$.  Suppose that $\psi \colon L_i^X \rightarrow L_j^Y$
; that $\phi \colon L_j^Y \rightarrow L_k^Z$; that and $\theta \colon L_i^X
\rightarrow L_k^Z$.  Then the \ct-morphisms 
\begin{align*}
\mu_{\psi} &\colon F(c,j,Y) \rightarrow F(b,i,X) \\
\mu_{\phi} &\colon F(d,k,Z) \rightarrow F(c,j,Y) \\
\mu_{\theta} &\colon F(d,k,Z) \rightarrow F(b,i,X) 
\end{align*}
defined as above satisfy $\mu_\theta = \mu_\psi \circ \mu_\phi$.
Hence if $\alpha$ is any element of $[F(b,i,X) \otimes\scrt M, N]_0$ then
$(\phi \circ \psi)(\alpha) = \theta(\alpha)$.
\end{proof}

\begin{prop} \label{adjoint}
If $L$, $M$, and $N$ are \ct-modules, then there is a natural isomorphism 
of graded groups
$$[L \otimes\scrt M,N]_* \cong [L, {\rm Hom}\scrt(M,N)]_*$$
\end{prop}

\begin{proof}
For simplicity, we will consider only homomorphisms of degree 0.  The 
isomorphism in arbitrary degree can be established analogously (or by 
suspending algebras).

We will define natural transformations
\begin{align*}
\Gamma &\colon [L \otimes\scrt M, N]_0 \rightarrow [L, \hom\scrt(M,N) ]_0 \\
\Delta &\colon [L, \hom\scrt(M,N)]_0 \rightarrow [L \otimes\scrt M, N]_0
\end{align*}
which are inverse to each other.

First we define $\Gamma$.  Let $\alpha \colon L \otimes\scrt M \rightarrow N$ and let
$l \in L_i^X$.  Then $\Gamma(\alpha)(l) \in \hom\scrt(M,N)_i^X = [F(b,i,X) \otimes
M, N]_0$ is defined to be the composite
$$F(b,i,X) \otimes M \xrightarrow{\mu_b^l \otimes 1}
	L \otimes M \xrightarrow{\alpha} N \; $$
where $\mu_b^l$ is the \ct-homomorphism which carries the generator $b \in F(b,i,X)$ to
the element $l \in L_i^X$.  The formula for $\Gamma$ is 
$\Gamma(\alpha)(l) = \alpha \circ (\mu_b^l \otimes 1)$.

To see that $\Gamma(\alpha)$ is a \ct-morphism, let $\theta$ be a \ct-operation sending $l
\in L_i^X$ to $\theta(l) \in L_j^Y$.  Let $b \in F(b,i,X)$ and let $c \in
F(c,j,Y)$.  Then  
\begin{align*}
\Gamma(\alpha)(\theta(l)) &= \alpha \circ (\mu_c^{\theta(l)} \otimes 1) \\
	&= \alpha \circ (\mu_b^l \otimes 1) \circ (\mu_c^{\theta(b)} \otimes 1)\\
	&= \alpha \circ (\mu_b^l \otimes 1) \circ (\mu_{\theta} \otimes 1) \\
	&= \Gamma(\alpha)(l) \circ (\mu_\theta \otimes 1) \\
	&= \theta(\Gamma(\alpha)(l))    \; . 
\end{align*}

Now we define $\Delta$.  Let $\beta \colon L \rightarrow \hom\scrt(M,N)$ be given. 
For any $l \in L_i^X$, the element $\beta(l)$ is in $\hom\scrt(M,N)_i^X =
[F(b,i,X) \otimes M, N]_0$.  So we can define $\Delta(\beta) \in [L \otimes M, 
N]_0$ by
the formula
  $$\Delta(\beta)(l \otimes m) = \beta(l)(b \otimes m) \;. $$
on the pure tensors of $L \otimes M$.  
To show that $\Delta(\beta)$ is an element of
$[L \otimes M, N]$, we must show that it is a
\ct-pairing.  Most of the \ct-pairing relations are of the form 
$$\theta_1(\Delta(\beta)(\theta_2(l) \otimes \theta_3(m)))
	= \psi_1 (\Delta(\beta)(\psi_2(l) \otimes \psi_3(m)))$$
for certain \ct-operations $\theta_i$ and $\psi_i$.  To show this holds, we use the
fact that $\beta(l)$ is a \ct-pairing for any $l \in L_i^X$.  Let $\theta_2(l) \in
L_j^Y$ and let $c$ be the generator of $F(c,j,Y)$.  Then 
\begin{align*}
\theta_1(\Delta(\beta)(\theta_2(l) \otimes \theta_3(m))) 
	&= \theta_1(\beta(\theta_2(l))(c \otimes \theta_3(m)))  \\
	&= \theta_1(\theta_2(\beta(l)(c \otimes \theta_3(m)))  \\
	&= \theta_1(\beta(l)(\mu_{\theta_2}(c) \otimes \theta_3(m)))  \\
	&= \theta_1(\beta(l)(\theta_2(b) \otimes \theta_3(m))) \\
	&= \psi_1(\beta(l)(\psi_2(b) \otimes \psi_3(m)))  \\
	&= \psi_1(\Delta(\beta)(\psi_2(l) \otimes \psi_3(m))) 
\end{align*}
The most general form of the \ct-pairing relations is
$$\theta_1 (\alpha(\theta_2(l) \otimes \theta_3(m)) )
	= \sum_{i = 1}^n (\psi_{i,1}(\alpha(\psi_{i,2}(l) \otimes \psi_{i,3}(m))))
\; $$  
which can be demonstrated similarly. 

Finally, we need to show that $\Gamma$ and $\Delta$ are inverses.
To show $\Gamma \circ \Delta = 1$, let $\beta \in [L, \hom\scrt(M,N)]_0$ be given. 
For any $l \in L_i^X$ and for any arbitrary pure tensor $\theta(b) \otimes m \in
F(b,i,X) \otimes M$ we have 
\begin{align*}
((\Gamma \circ \Delta)\beta)(l)(\theta(b) \otimes m)
	&= \Gamma(\Delta(\beta))(l)(\theta(b) \otimes m) \\
	&= (\Delta(\beta) \circ (\mu_b^l \otimes 1))(\theta(b) \otimes m) \\
	&= \Delta(\beta)(\theta(l) \otimes m)  \\
	&= \beta(\theta(l))(c \otimes m)  \\
	&= \theta(\beta(l))(c \otimes m)  \\
	&= \beta(l)(\mu_{\theta}(c) \otimes m)  \\
	&= \beta(l)(\theta(b) \otimes m)  \; .
\end{align*}

To show $\Delta \circ \Gamma = 1$, let $\alpha \in [L \otimes M, N]_0$ be given.  Then
\begin{align*}
(\Delta \circ \Gamma)(\alpha)(l \otimes m) 
	&= \Delta(\Gamma(\alpha))(l \otimes m)   \\
	&= \Gamma(\alpha)(l)(b \otimes m) \\
	&= \alpha \circ (\mu_b^l \otimes 1)(b \otimes m) \\
	&= \alpha(l \otimes m) \; .
\end{align*}
\end{proof}

For fixed a fixed \ct-module $M$ the co-variant functor $\hom\scrt(M,-)$ is left
exact because of its adjoint property (Theorem~2.6.1 in \cite{Weibel}).  
Therefore according to Definition~2.5.1 in \cite{Weibel}, and since by
Section~3.1 in \cite{Bou} there are enough injectives in \ct, we obtain the derived 
functor $\ext\scrt(M,N)$ which again is an element of \ct for any two \ct-modules
$M$ and $N$. 

The following propositions give us a more concrete representation of the 
real and complex parts of $\hom\scrt(M,N)$.  The isomorphisms given are of 
graded groups.

\begin{prop} \label{homo=[]}  If $M$ and $N$ are \ct-modules, then
$${\rm Hom}\scrt(M,N)^O \cong [M,N]_* \; .$$
\end{prop}

\begin{proof}
$$
\hom\scrt(M,N)_i^O 
= [F(b,i,\R) \otimes\scrt M, N]_0 
\cong [M, N]_i $$
\end{proof}

\begin{prop} \label{homu=homu} If $M$ and $N$ are \ct-modules, then
$${\rm Hom}\scrt(M,N)^U \cong {\rm Hom}_{KU_*(\R)}(M^U, N^U) \; .$$
\end{prop}

\begin{proof}
By use of suspensions, the problem is reduced to
that of showing that the two graded groups are isomorphic in graded degree zero.  By
definition,  $\hom\scrt(M,N)^U = [F(b,0,\C) \otimes M, N]_0$.  Recall from 
Proposition~3.6 of
\cite{Boer2} that   
\begin{align*}
F(b, 0, \C) \otimes\scrt M 
&\cong \{ \,
\{r(b \otimes m\su) \mid m\su \in M\su\},  \\
& \hspace{2cm} \{b \otimes m_1 + \psi\su (b \otimes m_2) \mid m_i \in M\su\},  \\
& \hspace{3cm} \{\gamma (b \otimes m_1) 
		+ \varepsilon r(b \otimes m_2) \mid m_i \in M\su \}  \}
\, \}    
\end{align*}
Any degree-0, \ct-module homomorphism from $F(b,0,\C) \otimes M$ to $N$ restricts to
a degree-0, $KU_*(\R)$-module homomorphism from $M$ to $N$ by restricting to
elements of the form $b \otimes m$ for $m \in M$.  This gives us a natural
transformation  
$$\Gamma \colon \hom\scrt(M,N)_0^U \rightarrow \hom_{KU_*(\R)}(M^U, N^U)_0 \; .$$

To define an inverse 
$$\Delta \colon \hom_{KU_*(\R)}(M^U, N^U)_0 \rightarrow \hom\scrt(M,N)_0^U \; $$
let $\phi \colon M^U \rightarrow N^U$ be a degree-0, $KU_*(\R)$-module homomorphism.
 Then define $\Gamma(\phi) \in [F(b,0,\C) \otimes M, N]$ by
\begin{align*}
\Delta(\phi) r(b \otimes m) &= r \phi(m) \\
\Delta(\phi)(b \otimes m) &= \phi(m)  \\
\Delta(\phi)\phi\su(b \otimes m) &= \psi\su \phi(m)  \\
\Delta(\phi)\gamma(b \otimes m) &= \gamma \phi(m)  \\
\Delta(\phi)\varepsilon r(b \otimes m) &= \varepsilon r \phi(m)
\end{align*}

It is clear that $\Gamma$ and $\Delta$ are inverses once we are convinced that
$\Delta(\phi)$ is an honest \ct-module homomorphism.  We must show that all of the
\ct-operations commute with $\Delta(\phi)$.  This is fairly straightforward but
involves many tedious computations.  As a representative example
we will show that the operation $\omega$ commutes with $\Delta(\phi)$, when 
operating on elements of the form $\gamma(b \otimes m)$ and $\ve r(b 
\otimes m)$. 
\begin{align*}
\Delta(\phi) \omega \gamma(b \otimes m)
	&= \Delta(\phi) \beta\st \gamma \zeta \gamma (b \otimes m) \\
	&= 0 \\
	&= \beta\st \gamma \zeta \gamma \phi(m) \\
	&= \omega \gamma \phi(m) \\
	&= \omega \Delta(\phi) \gamma(b \otimes m) \, 
\end{align*}
and
\begin{align*}
\Delta(\phi) \omega \ve r(b \otimes m) 
	&= \Delta(\phi) \beta\st \gamma \zeta \ve r(b \otimes m)  \\
	&= \Delta(\phi) \beta\st \gamma(1 + \psi\su)(b \otimes m) \\
	&= \Delta(\phi) 2 \gamma \beta\su^2 (b \otimes m)  \\
	&= \Delta(\phi) \gamma(b \otimes 2 \beta\su^2 m) \\
	&= \gamma \phi(2 \beta\su^2 m) \\
	&= 2 \beta\st \gamma \phi(m) \\
	&= \beta\st \gamma(1 + \psi\su) \phi(m) \\
	&= \beta\st \gamma \zeta \ve r \phi(m) \\
	&= \omega \Delta(\phi) \ve r (b \otimes m) \, .
\end{align*}
\end{proof}


\subsection{Hom and free objects} \label{homfree}

Recall that if $M$ is a $\Z$-graded object, then $\Sigma M$ is the 
$\Z$-graded object defined by $(\Sigma M)_n = M_{n+1}$.  Now, 
given a {\it CRT}-module $N$, let $\mathfrak{C}(N)$ be the {\it CRT}-module defined
to be the triple of graded groups, $\{N\pu, N\pu \oplus N\pu, \Sigma N\pu
\oplus N\pu \}$ with {\it CRT}-operations
\begin{align*}
\ve(n) &= (0,n) \\
\zeta(n_1, n_2) &= (n_2, n_2) \\
\psi\su(n_1, n_2) &= (n_2, n_1) \\
\psi\st(n_1, n_2) &= (-n_1, n_2) \\
\gamma(n_1, n_2) &= (n_1 + n_2, 0) \\
\tau(n_1, n_2) &= n_1
\end{align*}
with $c$ and $r$ defined by $c = \zeta \ve$ and $r = \tau \gamma$.
It is convenient to write the elements of this \ct-module formally as
\begin{align*}
\mathfrak{C}(N) 
&= \{ \,
\{ r(n\su) \mid n\su \in N\pu \}, \\
& \hspace{2cm} \{n\su \mid n\su \in N\pu \}
		\oplus \{\psi\su(n\su) \mid n\su \in N\pu \} \; , \\
& \hspace{3cm}   \{ \gamma(n\su)  \mid n\su \in N\pu \}
		 \oplus \{\ve r(n\su) \mid n\su \in N\pu \}  \, \} \; .
\end{align*}

Also, $\mathfrak{T}(N)$ is defined to be the triple $\{N\pt, N\pu
\oplus \Sigma N\pu, N\pt \oplus \Sigma N\pt \}$ with operations as given below:
\begin{align*}
\ve(n) &= (n, (-1)^{|n|} \ve \tau(n) ) \displaybreak[0]\\
\zeta(n_1, n_2) &= (\zeta(n_1), \zeta(n_2)) \displaybreak[0]\\
\psi\su(n_1, n_2) &= (\psi\su(n_1), \psi\su(n_2)) \displaybreak[0]\\
\psi\st(n_1, n_2) &=
	(\psi\st(n_1) + (-1)^{(|n_2|+1)} \beta\st^{-1} \omega(n_2),  
		\psi\st(n_2))\displaybreak[0] \\ 
\gamma(n_1, n_2) &= (\gamma(n_1), \gamma(n_2)) \displaybreak[0] \\ 
\tau(n_1, n_2) &= \ve \tau(n_1) + \eta\st(n_1) + (-1)^{{|n_2|}} n_2  
\end{align*}
In these formulas, the notation $|n|$ refers to the inherent graded degree of the
element $n \in N$.

\begin{prop} \label{homR}
 Let $N$ be an arbitrary {\it CRT}-module.  Then there is a natural 
 isomorphism
$$
{\rm Hom}\scrt(F(b, 0, \R) , N)  \cong N 
$$ 
\end{prop}

\begin{prop} \label{homC} 
 Let $N$ be an arbitrary {\it CRT}-module.  Then there is a natural 
 isomorphism
$$
{\rm Hom}\scrt(F(b, 0, \C) , N) \cong \mathfrak{C}(N)  
$$ 
\end{prop}

\begin{prop} \label{homT}
 Let $N$ be an arbitrary {\it CRT}-module.  Then there is a natural 
 isomorphism
$$
{\rm Hom}\scrt(F(b, 0, T) , N) \cong \mathfrak{T}(N) 
$$ 
\end{prop}

\begin{proof}[Proof of Proposition~\ref{homR}]
Let $F = F(b,0,\R)$.  For any $i \in \Z$ and $X \in \{O,U,T\}$ we have an 
isomorphism 
$$N_i^X \cong [F(c,i,X), N]_0 \cong [F(c,i,X) \otimes F, N]_0 = 
\hom\scrt(F,N)_i^X$$
given by $\Gamma \colon f \mapsto f(c \otimes b)$ 
(or simply $f \mapsto f(b)$).
So $N$ and $\hom\scrt(F,N)$ are isomorphic as graded groups.  To show
that the $\ct$-operations coincide, it suffices to show that $\Gamma$
commutes with a typical $\ct$-operation $\theta
\colon N_i^X \rightarrow N_j^Y$.  Given any $f \in [F(c,i,X), N]_0$:
$$\Gamma(\theta(f)) = (\theta(f))(c) = (f \circ \mu_\theta)(c)
	= f(\theta(c)) = \theta(f(c)) = \theta(\Gamma(f)) \; .$$
\end{proof}

\begin{proof}[Proof of Proposition~\ref{homC}]
Let $F = F(b,0,\C)$.  Then for any integer $i$ we have
\begin{align*}
\hom\scrt(F,N)_i^O &= [F(c,i,\R) \otimes F, N]_0 \\
		&= [F(x,i,\C), N]_0 \\
		&= N_i^U   \\
\hom\scrt(F,N)_i^U &= [F(c,i,\C) \otimes F, N]_0 \\
		&= [F(y,i,\C) \oplus F(z,i,\C), N]_0 \\
		&= N_i^U \oplus N_i^U  \\
\hom\scrt(F,N)_i^T &= [F(c,i,T) \otimes F, N]_0 \\
		&= [F(u,i,\C) \oplus F(v,i+1,\C), N]_0 \\
		&= N_i^U \oplus \Sigma N_{i}^U  \\
\end{align*}
since $F(c,i,\R) \otimes F(b,0,\C)$ has free generator $x = c(c) \otimes b$; the
\ct-module $F(c,i,\C) \otimes F(b,i,\C)$ has free generators $y = c \otimes b$ and
$z = \psi\su(c) \otimes b$; and $F(c,i,T) \otimes F(b,i,\C)$ has free generators $u =
\zeta(c) \otimes b$ and $v = c \tau(c) \otimes b$ (using Propositions~3.5, 3.6, and
3.7 of \cite{Boer2}).  Again, we have the isomorphism we want on the level of graded
groups.  We will show that the isomorphism respects the \ct-operations.

For any $n \in N\pu_i$, let $n^x \in [F(c,i,\R) \otimes F, N]_0 = \hom\scrt(F,N)\po_i$
be the \ct-homomorphism which carries $x$ to $n$.  In $[F(c,i,\C) \otimes F, N]_0 =
\hom\scrt(F,N)\pu_i$, let $n^y$ (respectively $n^z$) be the \ct-homomorphism which
carries $y$ (respectively $z$) to $n$ and $z$ (respectively $y$) to $0$.  For $n \in
N\pu_i$, let $n^u$ (respectively $n^v$) be the element of 
$[F(c,i,T) \otimes F, N]_0 = \hom\scrt(F,N)\pt_i$
which carries $u$ (respectively $v$) to $n$ and vanishes on $v$ 
(respectively $u$).  Finally, for $n \in N\pu_{i+1}$, let
$n^v \in [F(c,i,T) \otimes F, N] = \hom\scrt(F,N)\pt_i$ be the \ct-homomorphism
which carries $v$ to $n$ and vanishes on $u$.

The four relations
\begin{align*} r(n^y) &= n^x \\ 
		\gamma(n^y) &= n^v \\
		\ve(n^x) &= n^u \\
		\psi\su(n^y) &= n^z
\end{align*}
are demonstrated below and are sufficient to establish that the two
\ct-modules under consideration are isomorphic.  
\begin{align*}
r(n^y)(x) &= n^y \circ (\mu_r \otimes 1)(c(c) \otimes b) \\
	&= n^y(cr(c) \otimes b)  \\
	&= n^y((1 + \psi\su)(c) \otimes b) \\
	&= n^y(y + z) \\
	&= n  \displaybreak[0] \\
\gamma(n^y)(u) &= n^y(\mu_\gamma \otimes 1)(\zeta(c) \otimes b) \\
	&= n^y(\zeta(\gamma(c)) \otimes b) \\
	&= 0  \\
\gamma(n^y)(v) &= n^y(\mu_\gamma \otimes 1)(c \tau(c) \otimes b) \\
	&= n^y(c \tau \gamma(c) \otimes b) \\
	&= n^y((1 + \psi\su)(c) \otimes b) \\
	&= n^y(y + z) \\
	&= n \displaybreak[0] \\
\ve(n^x)(u) &= n^x(\mu_\ve \otimes 1)(\zeta(c) \otimes b) \\
	&= n^x(\zeta \ve (c) \otimes b) \\
	&= n^x(c(c) \otimes b) \\
	&= n  \\
\ve(n^x)(v) &= n^x(\mu_\ve \otimes 1)(c \tau(c) \otimes b) \\
	&= n^x(c \tau \ve(c) \otimes b) \\
	&= 0  \displaybreak[0]\\
\psi\su(n^y)(y) &= n^y(\mu_\psi \otimes 1)(c \otimes b) \\
	&= n^y(\psi\su(c) \otimes b) \\
	&= n^y(z) \\
	&= 0  \\
\psi\su(n^y)(z) &= n^y(\mu_\psi \otimes 1)(\psi\su(c) \otimes b) \\	
	&= n^y(y) \\
	&= n
\end{align*}
\end{proof}

\begin{proof}[Proof of Proposition~\ref{homT}]
Let $F = F(b,0,T)$.  Then we have
\begin{align*}
\hom\scrt(F,N)_i\po &= [F(c, i, \R) \otimes F, N] \\
	&= [F(x,i,T), N] \\
	&= N_i\pt  \\
\hom\scrt(F,N)_i\pu &= [F(c,i,\C) \otimes F, N]   \\
	&= [F(y,i,\C) \oplus F(z,i+1,\C), N]  \\
	&= N_i\pu \oplus \Sigma N_{i}\pu    \\
\hom\scrt(F,N)_i\pt &= [F(c,i,T) \otimes F, N]   \\
	&= [F(u,i,T) \oplus F(v, i+1, T), N]  \\
	&= N_i\pt \oplus \Sigma N_{i} \pt
\end{align*}
since, again referring to Section~3 of \cite{Boer2}, $F(c,i,\R) \otimes F(b,0,T)$
has a free generator $x = \ve(c) \otimes b$;  the \ct-module $F(c,i,\C) \otimes
F(b,0,T)$ has free generators $y = c \otimes \zeta b$ and $z = c \otimes c \tau b$;
and $F(c,i,T) \otimes F(b,0,T)$ has free generators $u = c \otimes b$ and $v = c
\otimes \ve \tau b$.  

As in the previous proof, we designate the following elements: $n^x \in \hom\scrt(F,N)\po_i$
for any $n \in N\pt_i$; $n^y \in \hom\scrt(F,N)\pu_i$ for $n \in
N\pu_i$; $n^z \in \hom\scrt(F,N)\pu_i$ for any $n \in N\pu_{i+1}$; $n^u \in \hom\scrt(F,N)\pt_i$ 
for any $n \in N\pt_i$; and finally 
$n^v \in \hom\scrt(F,N)\pt_i$ for any $n \in N\pt_{i+1}$.

To show that we have
an isomorphism of \ct-modules we need to show that the relations in the definition of
$\mathfrak{T}(N)$ hold for these elements defined in the previous paragraph.  
\begin{align*}
\ve(n^x) &= n^u + (-1)^i (\ve \tau(n))^v = n^u + (-1)^{|n|} (\ve \tau(n))^v \\ 
\zeta(n^u) &= (\zeta n)^y \\
\zeta(n^v) &= (\zeta n)^z \\
\psi\su(n^y) &= (\psi\su n)^y \\
\psi\su(n^z) &= (\psi\su n)^z \\
\psi\st(n^u) &= (\psi\st n)^u \\
\psi\st(n^v) &= (-1)^{i} \beta\st^{-1} \omega n^u + (\psi\st n)^v 
	= (-1)^{(|n|+1)} \beta\st^{-1} \omega n^u + (\psi\st n)^v     \\
\gamma(n^y) &= (\gamma n)^z \\
\gamma(n^z) &= (\gamma n)^z \\
\tau(n^u) &= (\ve \tau n + \eta\st n)^x \\
\tau(n^v) &= (-1)^{i+1} n^x = (-1)^{|n|} n^x
\end{align*}

We will work out a few of these, leaving the rest for the reader.  The calculations
\begin{align*}
\ve(n^x)(u) &= n^x(\ve(c) \otimes b) = n \\
\ve(n^x)(v) &= n^x(\ve(c) \otimes \ve \tau b) \\
	&= \ve n^x(c \otimes \tau b) \\
	&= (-1)^i \ve \tau n^x(\ve(c) \otimes b) \\
	&= (-1)^i \ve \tau(n) \displaybreak[0] \\
\tau(n^u)(x) &= n^u(\ve \tau(c) \otimes b) \\
	&= n^u(\ve \tau(c \otimes b) + (-1)^{i+1} c \otimes \ve \tau b 
					+ \eta\st c \otimes b) \\
	&= \ve \tau n + \eta\st n \\
\tau(n^v)(x) &= n^v(\ve \tau(c) \otimes b) \\
	&= n^v(\ve \tau(c \otimes b) + (-1)^{i+1} c \otimes \ve \tau b 
					+ \eta\st c \otimes b) \\
	&= (-1)^{i+1} n
\end{align*}
show that $\ve(n^x) = n^u + (-1)^i (\ve \tau n)^v$; $\tau(n^u) = (\ve \tau n +
\eta\st n)^x$; and $\tau(n^v) = (-1)^{i+1} n^x$.
\end{proof}


\subsection{Hom and acyclic objects} \label{homacyc}

\begin{prop} \label{freeacyclic}
Let $M$ be a free \ct-module and $N$ be an acyclic \ct-module.  
Then ${\rm Hom}\scrt(M,N)$
is acyclic. 
\end{prop}

\begin{proof}
First we consider the case that $M$ is monogenic.  If $M = F(b,0,\R)$, then by
Proposition~\ref{homR} $\hom\scrt(M,N) \cong N$ which is acyclic.

If $M = F(b,0,\C)$, then $\hom\scrt(M,N) \cong \mathfrak{C}(N)$ by
Proposition~\ref{homC}.  We claim that $\mathfrak{C}(N)$ is acyclic.  In fact, by
comparing to Proposition~3.6 in \cite{Boer2}, the \ct-module $\mathfrak{C}(N)$ is
isomorphic to $F(b,0,\C) \otimes N$, which is acyclic by Lemma~3.11 in \cite{Boer2}.

If $M = F(b,0,T)$, then $\hom\scrt(M,N) \cong \mathfrak{T}(N)$ by
Proposition~\ref{homT}.  By Proposition~3.7 in \cite{Boer2},
\begin{align*}
F(b, 0, T) \otimes\scrt N 
&\cong \{ \,
\{ \tau(b \otimes n\st) \mid n\st \in N\st \},  \\
& \hspace{2cm} \{ \zeta b \otimes n\su \mid n\su \in N\su \} \oplus
	\{ c \tau b \otimes n\su \mid n\su \in N\su \},  \\
& \hspace{3cm} \{ b \otimes n\st \mid n\st \in N\st \} \oplus 
	\{ \varepsilon \tau b \otimes n\st \mid n\st \in N\st \} \, \} \\
&\cong \{\Sigma^{-1} N\st, N\su \oplus \Sigma^{-1} N\su, N\st \oplus 
	\Sigma^{-1} N\st \}   \; . 
\end{align*}
It can easily be shown that the defining relations of $\mathfrak{T}(N)$ are
satisfied by the elements 
\begin{align*}
(-1)^{|n\st|} \tau(b \otimes n\st) &\in (F(b, 0, T) \otimes\scrt N)\po \\
(-1)^{|n\su|} c \tau b \otimes n\su &\in (F(b, 0, T) \otimes\scrt N)\pu \\
\zeta b \otimes n\su &\in (F(b, 0, T) \otimes\scrt N)\pu  \\
(-1)^{|n\st|} (\ve \tau b \otimes n\st + \eta\st(b \otimes n\st)) 
		&\in (F(b, 0, T) \otimes\scrt N)\pt \\
b \otimes n\st &\in (F(b, 0, T) \otimes\scrt N)\pt
\end{align*}
of $F(b,0,T) \otimes\scrt N$.  Therefore, $\mathfrak{T}(N)$ is isomorphic to
$F(b,0,T) \otimes\scrt N$ (with a shift); and so Lemma~3.13 in 
\cite{Boer2} implies that
$\mathfrak{T}(N)$ is acyclic.

Now the suspension formula
$$\hom\scrt(\Sigma M, N) = \Sigma^{-1} \hom\scrt(M,N)$$
follows directly from the definition and implies that $\hom\scrt(M,N)$ is acyclic if
we take $M$ to be any free monogenic \ct-module.  

Finally, we can write an arbitrary free \ct-module as a direct sum of monogenic
free \ct-modules (see Section~2.4 of \cite{Bou}).  Let $M = \bigoplus_{\lambda \in
\Lambda} M_\lambda$.  Then we claim that $\hom\scrt(M,N) \cong \prod_{\lambda \in
\Lambda} \hom\scrt(M_\lambda, N)$.  It then follows that $\hom\scrt(M,N)$ is acyclic.
To prove the claim, we have
\begin{alignat*}{2}
\hom\scrt \left( \bigoplus_{\lambda \in \Lambda} M_\lambda, N \right)^X_i
	&\cong \left[ F(b,i,X) \otimes \bigoplus_{\lambda \in \Lambda} M_\lambda, 
			N \right]  \\ 	
	&\cong \left[ \bigoplus_{\lambda \in \Lambda} 
			\left(F(b,i,X) \otimes M_\lambda \right), N \right]  
	 	& & \qquad \text{by Lemma~IV.1.6 in \cite{Boer} } \\
	&\cong \prod_{\lambda \in \Lambda} \left[ F(b,i,X) \otimes M_\lambda, 
			N \right]  		
		& & \qquad \text{by Exercise~A.1.4 in \cite{Weibel}} \\ 	
	&\cong \prod_{\lambda \in \Lambda} \hom\scrt(M_\lambda,
			N)_i^X \; . 
\end{alignat*}
Since these isomorphisms are natural, they commute with respect to any homomorphism
$\mu_\theta \colon F(c,j,Y) \rightarrow F(b,i,X)$ and so we have a \ct-module
homomorphism as desired.
\end{proof}

\begin{prop} \label{injacyclic}
Let $M$ be an acyclic \ct-module and let $N$ be an injective \ct-module.  Then
${\rm Hom}\scrt(M,N)$ is acyclic. 
\end{prop}

\begin{proof}
Because $M$ is acyclic, it has projective dimension at most 1 (Theorem~3.4 in
\cite{Bou}) so we can form a resolution of $M$
$$0 \rightarrow F_1 \rightarrow F_0 \rightarrow M \rightarrow 0$$
where $F_0$ and $F_1$ are projective.  Then by Theorem~3.2 in \cite{Bou} $F_0$ and
$F_1$ are free.  Now, applying $\hom\scrt(-,N)$ to this resolution, we obtain
$$\ext\scrt(M,N) \leftarrow \hom\scrt(F_1,N) \leftarrow \hom\scrt(F_0,N) 
	\leftarrow \hom\scrt(M,N) \leftarrow 0 \; .$$
Since $N$ is injective, $\ext\scrt(M,N) = 0$ so we have a short exact sequence of
\ct-modules.

By Theorem~3.3 in \cite{Bou}, the injectivity of $N$ implies it is
acyclic.  Then by Proposition~\ref{freeacyclic},
$\hom\scrt(F_0,N)$ and $\hom\scrt(F_1,N)$ are acyclic.  Finally, the Snake Lemma
(1.3.2 in \cite{Weibel}) implies $\hom\scrt(M,N)$ is acyclic.
\end{proof}


\section{The Universal Coefficient Theorem}

By Proposition~\ref{adjoint}, the \ct-pairing
$$\alpha \colon 
KK\crt(\R, A) \otimes\scrt KK\crt(A, B) \longrightarrow 
		KK\crt(\R, B)$$
described in Section~\ref{intproduct}
implies the existence
of a natural \ct-homomorphism
$$\gamma \colon KK\crt(A,B) \rightarrow \hom\scrt(K\crt(A), K\crt(B))$$
which is the adjoint of $\alpha$.  We note that we are also implicitly 
making use of the isomorphism between $KK\crt(\R, A) \otimes\scrt KK\crt(A, B)$ 
and $KK\crt(A,B) \otimes\scrt KK\crt(\R,A)$.

\begin{thm} \label{mainuct}
Let $A$ and $B$ be real separable C*-algebras 
such that $\C \otimes A$ is in the bootstrap category $\mathcal{N}$.  
Then there is a short
exact sequence
\end{thm}
$0 \rightarrow {\rm Ext}\scrt(K\crt(A), K\crt(B)) \xrightarrow{\kappa}
	KK\crt(A,B) \xrightarrow{\gamma}
	{\rm Hom}\scrt(K\crt(A), K\crt(B)) 
	\rightarrow 0$
where $\kappa$ has degree $-1$.

As in Rosenberg and Schochet in \cite{RS}, the proof of this theorem will 
proceed in two steps.  The first step (in Section~\ref{specialcasesection})
addresses the special case in which 
$K\crt(B)$ is injective and the second step (in 
Section~\ref{generalcasesection}) 
uses a geometric injective 
resolution of $K\crt(B)$ to extend the result to arbitrary separable $B$.  
Our debt to 
Rosenberg and Schochet goes beyond analogy, however, since we 
will actually reduce our proof in the special case to an application 
of the universal coefficient theorem for complex C*-algebras.  

In section~\ref{corsection}, 
we will discuss the most immediate ramifications of 
Theorem~\ref{mainuct}.  Throughout the rest of this paper, all C*-algebras 
are assumed to be separable.

\subsection{Proof of the UCT --- Special Case} \label{specialcasesection}

In the special case that $K\crt(B)$ is injective, the \ct-module 
$\hom\scrt(K\crt(A), K\crt(B))$ vanishes.  Therefore the universal 
coefficient theorem collapses to the following.

\begin{prop} \label{specialcase}
Let $A$ and $B$ be real C*-algebras such that $\C \otimes A$ is in 
$\mathcal{N}$ and $K\crt(B)$ is an injective \ct-module. 
Then
$$\gamma \colon KK\crt(A,B) \rightarrow {\rm Hom}\scrt(K\crt(A), K\crt(B))$$
is an isomorphism.
\end{prop}

Before we prove this proposition, we need two lemmas relating $KK$-theory
for real C*-algebras to $KK$-theory for complex C*-algebras.  Up to this
point, we have been dealing entirely with real $KK$-theory which for real
C*-algebras $A$ and $B$ is defined in terms of real Kasparov bimodules as
in \cite{Schroder}, Section 2.3.  Recall that a real Kasparov
$(A,B)$-bimodule is a triple $(E, \phi, T)$ consisting of a Hilbert
$B$-module $E$, a graded *-homomorphism $\phi \colon A
\rightarrow \mathcal{L}_B(E)$, and a degree one operator $T \in 
\mathcal{L}_B(E)$ such that the elements
$(T - T_*)\phi(a)$, $(T^2 - 1)\phi(a)$, and $[T, \phi(a)]$ all lie in
$\mathcal{K}_B(E)$ for all $a \in A$.  In particular, even if $A$ and $B$
happen to be complex C*-algebras, the homomorphism $\phi$ is only required
to be real.  Note that if $B$ is a complex C*-algebra, then $E$ 
(and hence $\mathcal{L}_B(E)$) inherits a complex structure using the 
approximate identity of $B$ together with the fact that the approximate 
identity for $B$ is an approximate identity for the action of $B$ on $E$ 
(Lemma 1.1.4 of \cite{JT}).

By contrast, if $A$ and $B$ are complex C*-algebras, let $KK\sc(A,B)$
denote the complex $KK$-theory comprised of equivalence classes of complex
Kasparov bimodules $(E, \phi, T)$.  This is the $KK$-theory of Rosenberg
and Schochet's Universal Coefficient Theorem.  The difference is that now
$\phi$ is required to respect the complex structures of $A$ and
$\mathcal{L}_B(E)$.  Since every complex homomorphism is a real
homomorphism, there is a natural transformation $\rho \colon KK\sc(A,B)
\rightarrow KK(A,B)$ which forgets the complex structures.  The next lemma 
describes another natural transformation connecting the two versions of 
$KK$-theory.

\begin{lemma} \label{kku}
Let $A$ be a real C*-algebra and let $B$ be a complex C*-algebra.  There 
is a natural isomorphism 
$\nu \colon KK\sc(\C \otimes A,B) \rightarrow KK(A, B)$.
\end{lemma}

\begin{proof}
Let $(E, \phi, T)$ be a complex Kasparov $(\C \otimes A, B)$-bimodule
representing an element of $KK\sc(\C \otimes A,B)$.  Then by restricting
$\phi$ to $\phi|_A \colon A \rightarrow \mathcal{L}_B(E)$, we obtain a real
Kasparov $(A, B)$-bimodule $(E, \phi|_A, T)$.  This mapping induces a
homomorphism $\nu$ which has an inverse: given a real Kasparov 
$(A,B)$-bimodule $(E, \phi, T)$, the homomorphism 
$\phi \colon A \rightarrow \mathcal{L}_B(E)$ 
extends uniquely to a complex homomorphism 
$\phi\sc \colon \C \otimes A \rightarrow \mathcal{L}_B(E)$ 
using the complex structure of $\mathcal{L}_B(E)$.
\end{proof}

Now, we recall the mechanics of the intersection product described for the 
general case in \cite{Kasp}.  Other expositions include Section~2.4 in 
\cite{Schroder} for real C*-algebras; and Section~18.3 in \cite{Black} and 
Section~2.2 in \cite{JT} for complex C*-algebras.  Let $A$, $B$, and $C$ 
be real C*-algebras; let $x \in KK(A,B)$ be represented by an 
$(A,B)$-bimodule $(E, \phi, T)$; and let $y \in KK(B,C)$ be represented by 
a $(B,C)$-bimodule $(F, \psi, S)$.  Then the product $\alpha(x \otimes 
y) \in KK(A,C)$ is given by an $(A,C)$-bimodule $(E \otimes_\psi F, \phi 
\otimes_\psi 1, T \# S)$.  The operator $T \# S$ can be taken to be any 
operator $U$ which satisfies:
\begin{align} 
T_x S - (-1)^{\partial x} UT_x &=0 && \mod \mathcal{K}(F, E \otimes_\psi F) 
		\notag \\
S T_x^* - (-1)^{\partial x} T_x^*U &=0 && \mod \mathcal{K}(E \otimes_\psi F , E) 
		\label{connection}  \\
\phi(a) [T \otimes 1, U] \phi(a)^* &\geq 0 
	\qquad && \mod \mathcal{K}(E \otimes_\psi F)    \notag
\end{align}
for all $x \in E$ and $a \in A$.  Here $T_x \in \mathcal{L}(F, E \otimes_\psi F)$ 
defined by $T_x(y) = x \otimes y$ and $\partial x$ is the graded degree 
of $x$.
This defines a pairing
$$\alpha \colon KK(A,B) \otimes KK(B,C) \rightarrow KK(A,C) \; .$$
The same construction also defines a pairing
$$\alpha\sc \colon KK\sc(A,B) \otimes KK\sc(B,C) \rightarrow KK\sc(A,C) $$
if $A$, $B$, and $C$ are complex C*-algebras (not to be confused with 
$\alpha\su$ described in Section~\ref{intproduct}).

We can also effect a pairing
$$\alpha'
\colon KK(A, \C \otimes B) \otimes KK(B, \C \otimes C)
\rightarrow KK(A, \C \otimes C)$$
by $\alpha' = \alpha \circ (1 \otimes \mu_*) \circ (1 \otimes \tau\sc)$
where $\mu \colon \C \otimes \C \otimes C \rightarrow \C \otimes C$ is 
given by complex multiplication and the $\tau\sc$ is the homomorphism
$$KK(B, \C \otimes C) \rightarrow KK(\C \otimes B, \C \otimes \C \otimes 
C)$$ defined by $(E, \phi, T) \mapsto (\C \otimes E, 1 \otimes \phi, 1 
\otimes T)$ (according to the external tensor product construction in 
Section~2.1.6 of \cite{JT}).

\begin{lemma} \label{nucommutes}
Let $A$, $B$, and $C$ be real C*-algebras.  
Then the following diagram commutes:
\end{lemma}
$$ \xymatrix{
KK\sc(\C \otimes A, \C \otimes B) \otimes KK\sc(\C \otimes B, \C \otimes C) 
			\ar[rr]^-{\alpha\sc} \ar[d]^{\nu \otimes \nu}	&&
KK\sc(\C \otimes A, \C \otimes C) 
			\ar[d]^{\nu}					\\
KK(A, \C \otimes B) \otimes KK(B, \C \otimes C) 
			\ar[rr]^-{\alpha'}			&&
KK(A, \C \otimes C)
									}$$

\begin{proof}
We claim first that $\mu_* \circ \tau_\C \circ \nu = \rho
\colon KK\sc(\C \otimes B, \C \otimes C \rightarrow KK(\C \otimes B, \C 
\otimes C)$.  Let $x$ be an element 
of $KK\sc(\C \otimes B, \C \otimes C)$ represented by a complex Kasparov 
$(\C \otimes B, \C \otimes C)$-bimodule $(F, \psi, S)$.  Then $\nu(x)$ is 
represented by $(F, \psi|_B, S)$ and $(\tau_\C \circ \nu)(x) \in 
KK(\C \otimes B, \C 
\otimes \C \otimes C)$ is represented by $(\C \otimes F, 1 \otimes \psi|_B, 
1 \otimes S)$.  Finally  
$\mu_* \tau \nu(x)$ is given by 
$(\widetilde{F}, \mu_*(1 \otimes \psi_B), \mu_*(1 \otimes S))$ 
where $\widetilde{F}$ is the Hilbert $\C \otimes C$-module 
$$\widetilde{F} = \mu_*(\C \otimes F) 
= (\C \otimes F)/ \{x | \mu\langle x,x \rangle = 0 \} \; $$
(according to the push-out construction of Section~2.1.5 of \cite{JT}).
There is a Hilbert $\C \otimes C$-module homomorphism $\mu_F \colon \widetilde {F} 
\rightarrow F$ defined by $\mu_F [\lambda \otimes f] = \lambda f$.  
Since $\langle \mu_F(x), \mu_F(y) \rangle = \mu\langle x,y \rangle$, this is 
well-defined.  Define $\mu_F^{-1}(f) = [1 \otimes f]$.  Clearly $\mu_F 
\circ \mu_F^{-1} = 1$.  To show that $\mu_F^{-1} \circ \mu_F = 1$, it is 
straightforward to calculate that $\mu \langle x, x \rangle = 0$ if 
$x = \lambda \otimes f - 1 \otimes \lambda f$. 
Therefore, $\mu_F$ is an isomorphism.  Under this isomorphism, the 
homomorphism
$\mu_*(1 \otimes \psi_B)$ corresponds to $\psi$ and the operator 
$\mu_*(1 \otimes S)$ corresponds to $S$.  This proves the claim.  

Thus the diagram that we 
wish to commute can be redrawn as
$$ \xymatrix{
KK\sc(\C \otimes A, \C \otimes B) \otimes KK\sc(\C \otimes B, \C \otimes C) 
			\ar[rr]^-{\alpha\sc} \ar[d]^{\nu \otimes \rho}	&&
KK\sc(\C \otimes A, \C \otimes C) 
			\ar[d]^{\nu}					\\
KK(A, \C \otimes B) \otimes KK(\C \otimes B, \C \otimes C) 
			\ar[rr]^-{\widetilde{\alpha}}				&&
KK(A, \C \otimes C)			} \; .$$
To prove that this commutes, let $x$ and $y$ be elements of $KK\sc(\C \otimes A, 
\C \otimes B)$ and $KK\sc(\C \otimes B, \C \otimes C)$ respectively.  We 
must show that
$\alpha(\nu(x) \otimes \rho(y)) = \nu \alpha\sc(x \otimes y)$.  Let $x$
be represented by a Kasparov $(\C \otimes A, \C \otimes B)$-bimodule $(E, 
\phi, T)$ as before and let $y$ be represented by a Kasparov $(\C \otimes B, \C \otimes 
C)$-bimodule $(F, \psi, S)$.

Now, $\nu(x) = (E, \phi_A, T)$ and $\rho(y) = (F, \psi, S)$ are both real 
Kasparov bimodules, as is the product $\alpha(\nu(x) \otimes 
\rho(y)) = (E \otimes_\psi F, \phi_A \otimes 1, T \# S)$.  On the other 
hand, $x = (E, \phi, T)$ and $y = (F, \psi, S)$ are both complex Kasparov 
bimodules, as is the product $\alpha\sc(x \otimes y) = (E \otimes_\psi F, 
\phi \otimes 1, T \#S)$.  However, since $\psi$ respects the complex 
structures of $\C \otimes B$ and $\mathcal{L}_{\C \otimes C}(\C \otimes 
E)$, the Hilbert $(\C \otimes C)$-module $E \otimes_\psi F$ is the same in 
either case.  Furthermore, the operator $T \# S \in \mathcal{L}_{\C \otimes C}(E 
\otimes_\psi F)$ can be taken to be the same in either case since any
operator which satisfies the formulae given in (\ref{connection}) for all $a \in \C 
\otimes A$ will certainly satisfy the same formulae for all $a \in A$.
Finally, since the 
restriction of $\phi \otimes 1 \colon \C \otimes A \rightarrow 
\mathcal{L}_{\C \otimes C}(E \otimes_\psi F)$ to $A$ is $\phi|_A \otimes 
1$, it follows that $\nu(\alpha\sc(x \otimes y)) = \alpha(\nu(x) 
\otimes \rho(y))$.
\end{proof}

\begin{proof}[Proof of Proposition~\ref{specialcase}]
According to Propositions~\ref{KKcrt} and \ref{injacyclic},
both $KK\crt(A,B)$ and $\hom\scrt(K\crt(A), K\crt(B))$ are acyclic \ct-modules.
Hence by Section~2.3 of \cite{Bou}, to show that
$$\gamma \colon KK\crt(A,B) \rightarrow 
	\hom\scrt(K\crt(A), K\crt(B))$$ 
is an isomorphism, 
it suffices to show that the complex part 
$$KK\crt(A,B)\pu \rightarrow 
	\hom\scrt(K\crt(A), K\crt(B))\pu$$
is an isomorphism.  Using the definition of united $KK$-theory and 
Proposition~\ref{homu=homu}, this is equivalent to showing that the 
homomorphism
$$\gamma\su \colon KK(A, \C \otimes B) \rightarrow
		\hom_{KU_*(\R)}(K(\C \otimes A), K(\C \otimes B))$$
(which is the adjoint of $\alpha\su$) is an isomorphism.  

The following diagram commutes by Lemma~\ref{nucommutes}.
$$ \xymatrix{
KK\sc(\C \otimes A, \C \otimes B) \ar[rrr]^-{\gamma\sc} \ar[d]^{\nu}  &&& 
\hom_{KU_*(\R)}(KK\sc(\C, \C \otimes A), KK\sc(\C, \C \otimes B)) 
					\ar[d]^{\hom(\nu,\nu)} \\
KK(A, \C \otimes B) \ar[rrr]^-{\gamma\su}		  &&&
\hom_{KU_*(\R)}(K(\C \otimes A), K(\C \otimes B))
			}$$  
Since $K\crt(B)$ is an injective 
\ct-module, the complex part $K(\C \otimes B) = KK\sc(\C \otimes \C 
\otimes B$ consists of divisible groups (Theorem~3.3 
in \cite{Bou}).  Hence, Theorem~2.1 of 
\cite{RS} implies that $\gamma\sc$ is an isomorphism and therefore so is 
$\gamma\su$.

\end{proof}   

\subsection{Proof of the UCT --- General Case} \label{generalcasesection}

\begin{lemma} \label{Nexists}
There exists a real C*-algebra $N$ containing a projection $p \in N$ such that 
1) $K\crt(N)$ is an injective \ct-module and 
2) for all real C*-algebras such 
that $K\crt(F)$ is free, the homomorphism $F \rightarrow 
F \otimes N$ defined by $a \mapsto a \otimes p$ induces a monomorphism of 
united $K$-theory.
\end{lemma}

The C*-algebra $N$ will play a role analogous to the 
role in Theorem~3.2 of \cite{RS} played by the complex UHF-algebra whose 
$K$-theory is $\Q$ in even degrees.  In Table~\ref{rational}, we have the 
united $K$-theory of the real UHF-algebra $N_1$ with $K_0(N_1) \cong \Q$.  
However, because of the presence of 2-torsion in free \ct-modules, this 
algebra will not suffice.  To capture this 2-torsion, we will also involve 
a real C*-algebra $N_2$ whose united $K$-theory is given in Table~\ref{2div}.

\begin{table} 
\caption{$K\crt(N_1)$} \label{rational}
$$\begin{array}{|c|c|c|c|c|c|c|c|c|c|}  
\hline \hline  
n & 0 & 1 & 2 & 3 & 4 & 5 & 6 & 7 & 8 \\
\hline  \hline
KO_*(N_1)
& \Q & 0 & 0 & 0 
	& \Q & 0 & 0 & 0 & \Q \\
\hline  
KU_*(N_1)
& \Q & 0 & \Q  & 0 & \Q & 0 
	& \Q & 0 & \Q  \\
\hline  
KT_*(N_1) 
& \Q & 0 & 0 & \Q  
& \Q & 0 & 0 & \Q & \Q    \\
\hline \hline
c_n & 1 & 0 & 0 & 0 & 2 & 0 & 0 & 0 & 1    \\
\hline
r_n & 2 & 0 & 0 & 0 & 1 & 0 & 0 & 0 & 2      \\
\hline
\varepsilon_n & 1 & 0 & 0 & 0 & 2 & 0 & 0 & 0 & 1   \\
\hline
\zeta_n &  1 & 0 & 0 & 0 &  1 & 0 & 0 & 0 & 1      \\
\hline
(\psi\su)_n & 1 & 0 & -1 & 0 & 1 & 0 & -1 & 0 & 1    \\
\hline
(\psi\st)_n & 1 & 0 & 0 & -1 & 1 & 0 & 0 & -1 & 1   \\
\hline
\gamma_n & 1 & 0 & 0 & 0 & 1 & 0 & 0 & 0 & 1     \\
\hline
\tau_n & 0 & 0 & 0 & 1 & 0 & 0 & 0 & 2 & 0        \\
\hline \hline
\end{array}$$
\end{table}

\begin{table} 
\caption{$K\crt(N_2)$} \label{2div}
$$\begin{array}{|c|c|c|c|c|c|c|c|c|c|}  
\hline  \hline 
n & \makebox[1cm][c]{0} & \makebox[1cm][c]{1} & 
\makebox[1cm][c]{2} & \makebox[1cm][c]{3} 
& \makebox[1cm][c]{4} & \makebox[1cm][c]{5} 
& \makebox[1cm][c]{6} & \makebox[1cm][c]{7} 
& \makebox[1cm][c]{8} \\
\hline  \hline 
KO_*(N_2)  &
\Z(2^\infty) & 0 & \Z_2 & \Z_2 & 
\Z(2^\infty) & 0 & 0 & 0 & 
\Z(2^\infty)         \\
\hline  
KU_*(N_2)   &
\Z(2^\infty) & 0 & 
\Z(2^\infty) & 0 & 
\Z(2^\infty)& 0 & 
\Z(2^\infty) & 0 & 
\Z(2^\infty)         \\
\hline  
KT_*(N_2)   &
\Z(2^\infty) & 0 & \Z_2 & \Z(2^\infty) &  
\Z(2^\infty) & 0 & \Z_2 & \Z(2^\infty) &
\Z(2^\infty)         \\
\hline
\hline
c_n & 1 & 0 & \tfrac{1}{2} & 0 & 2 & 0 & 0 & 0 & 1    \\
\hline
r_n & 2 & 0 & 0 & 0 & 1 & 0 & 0 & 0 & 2      \\
\hline
\varepsilon_n & 1 & 0 & 1 & \tfrac{1}{2} & 2 & 0 & 0 & 0 & 1   \\
\hline
\zeta_n &  1 & 0 & \tfrac{1}{2} & 0 &  1 & 0 & \tfrac{1}{2} & 0 & 1      \\
\hline
(\psi\su)_n & 1 & 0 & -1 & 0 & 1 & 0 & -1 & 0 & 1    \\
\hline
(\psi\st)_n & 1 & 0 & 1 & -1 & 1 & 0 & 1 & -1 & 1   \\
\hline
\gamma_n & 1 & 0 & 0 & 0 & 1 & 0 & 0 & 0 & 1     \\
\hline
\tau_n & 0 & 0 & 1 & 1 & 0 & 0 & 0 & 2 & 0        \\
\hline \hline
\end{array}$$
\end{table}

\begin{proof}
We will first construct the algebra $N_2$ as a limit of a directed system of matrix
algebras over real Cuntz algebras.  We begin by proving the claim that for any
positive integers $n$ and $m$ there is a unital homomorphism 
$$\phi \colon \mathcal{O}_{n+1} \rightarrow M_m(\mathcal{O}_{nm+1})\; $$
using a construction similar to that in the proof of Proposition~2.5 of 
\cite{PS}.   
Let $\mathcal{O}_{n+1}$ be generated by the isometries $S_0, S_1,
\dots, S_n$ with orthogonal ranges and which satisfy $S_i^* S_i = 1$ and 
$\sum_{i=0}^n S_i S_i^* = 1$. 
Let $\mathcal{O}_{mn+1}$ similarly be generated by the isometries $T_0, T_{ij}$
for $1 \leq i \leq n$, $1 \leq j \leq m$ which satisfy similar relations. 
Now, we define elements $S'_0, S'_1, \dots, S'_n \in M_m(\mathcal{O}_{mn+1})$
by
$$
S'_0 = \left( \begin{matrix}  
		T_0 & 0 & 0 & \dots & 0 \\
		0  &  1 & 0 & \dots & 0 \\
		0  &  0 & 1 & \dots & 0 \\
		\vdots & \vdots & \vdots & \ddots & \vdots   \\
		0 & 0 & 0 & \dots & 1  
	\end{matrix}  \right)    \qquad
S'_i = \left( \begin{matrix}
		T_{i1} & T_{i2} & \dots & T_{im} \\
		0 & 0 & \dots & 0 \\
		\vdots & \vdots & \ddots & \vdots \\
		0 & 0 & \dots & 0
	\end{matrix} \right) \; \text{~for $i \geq 1$.}
$$
It can be checked easily that these elements satisfy $(S'_i)^* S'_i = 1$ and 
$\sum_{i=0}^n S'_i (S'_i)^* = 1$ and so there is a C*-algebra homomorphism
$\phi$ which takes each $S_i$ to $S'_i$ .  The homomorphism $\phi$ is unital
since $\phi(1) = \phi(S_1^* S_1) = (S'_1)^* S'_1 = 1$.  This proves the claim.

Then we define the algebra $N_2$ to be the limit of
the system $$       \mathcal{O}_{4+1} \rightarrow 
		M_2(\mathcal{O}_{8+1}) \rightarrow 
		M_4(\mathcal{O}_{16+1}) \rightarrow 
		M_8(\mathcal{O}_{32+1}) \rightarrow \dots \; $$
where each connecting map is $M_{2^{k-2}}(\phi)$ with $n = 2^{k}$ and $m=2$.  
The united $K$-theory of $N_2$ is the limit of the system of \ct-modules
$\{K\crt(\mathcal{O}_{2^k+1}), M_{2^{k-2}}(\phi)_*\}$.  The united $K$-theory of
$\mathcal{O}_{2^k+1}$ is given in Table~\ref{cuntzC}. 
Since $\phi$ induces multiplication by 2 on $KO_0(\mathcal{O}_{2^k+1})$, 
we obtain $KO_0(N_2) = \Z(2^{\infty})$, as in all the gradings 
of $K\crt(N_2)$ in which $K\crt(\mathcal{O}_{2^k+1}) \cong \Z_{2^k}$.

Now $\phi_* = 0$ on $KO_1(\mathcal{O}_{2^k+1}) \cong \Z_2$ since the 
nonzero element of that group is $\eta\so(1)$.  The same is true for 
$KT_2(\mathcal{O}_{2^k+1})$ (with nonzero element $\ve \eta\so(1)$),  the first
summand of $KO_2(\mathcal{O}_{2^k+1})$ ($\eta\so^2(1)$), and 
$KT_5(\mathcal{O}_{2^k+1})$ ($\beta\st \ve
\eta\so(1)$).  Therefore, these groups vanish
when we pass to the limit.  

\begin{table}   
\caption{$K\crt(\mathcal{O}_{2^k+1})$ for $k \geq 2$} \label{cuntzC}
$$\begin{array}{|c|c|c|c|c|c|c|c|c|c|}  
\hline  \hline  
n & \makebox[1cm][c]{0} & \makebox[1cm][c]{1} & 
\makebox[1cm][c]{2} & \makebox[1cm][c]{3} 
& \makebox[1cm][c]{4} & \makebox[1cm][c]{5} 
& \makebox[1cm][c]{6} & \makebox[1cm][c]{7} 
& \makebox[1cm][c]{8} \\
\hline  \hline 
KO_n  &
\Z_{2^k} & \Z_2 & \Z_2^2 & \Z_2 & \Z_{2^k} &
0 & 0 & 0 & \Z_{2^k}         \\
\hline  
KU_n   &
\Z_{2^k} & 0 & 
\Z_{2^k} & 0 & 
\Z_{2^k} & 0 & 
\Z_{2^k} & 0 & 
\Z_{2^k}         \\
\hline  
KT_n   &
\Z_{2^k} & \Z_2 & \Z_2 & \Z_{2^k} &  
\Z_{2^k} & \Z_2 & \Z_2 & \Z_{2^k} &
\Z_{2^k}         \\
\hline \hline
c_n & 1 & 0 & \smh{0}{2^{k-1}} & 0 & 2 & 0 & 0 & 0 & 1 \\
\hline
r_n & 2 & 0 & \smv{1}{0} & 0 & 1 & 0 & 0 & 0 & 2  \\
\hline
\varepsilon_n & 1 & 1 & \smh{0}{1} & 2^{k-1} & 2 & 0 & 0 & 0 & 1 \\
\hline
\zeta_n & 1 & 0 & 2^{k-1} & 0 & 1 & 0 & 2^{k-1} & 0 & 1 \\
\hline
(\psi\su)_n & 1 & 0 & -1 & 0 & 1 & 0 & -1 & 0 & 1 \\
\hline
(\psi\st)_n & 1 & 1 & 1 & -1 & 1 & 1 & 1 & -1 & 1 \\
\hline
\gamma_n & 1 & 0 & 1 & 0 & 1 & 0 & 1 & 0 & 1 \\
\hline
\tau_n & 1 & \smv{1}{0} & 1 & 1 & 0 & 0 & 0 & 2 & 1  \\
\hline \hline
\end{array}$$
\end{table}

On the second summand of
$KO_2(\mathcal{O}_{2^k+1})$, we claim that $\phi_*$ is an isomorphism.  
Indeed, let $\smv{0}{1}$ be the non-zero element
of  $KO_2(\mathcal O_{2^k+1})$ representing the
second summand and let $y$ be a generator of $KU_2(\mathcal{O}_{2^k+1}) =
\Z_k$.   Then 
$(c \circ \phi_*)\smv{0}{1} = (\phi_* \circ c)\smv{0}{1} 
= \phi_*(\tfrac{k}{2} \cdot y) 
	= k \cdot y$.  
Therefore $\phi_*\smv{0}{1} = \smv{\alpha}{1}$ where $\alpha = 0$ or $1$. 
Using the naturality of the operations $\ve$ and $\eta\so$, 
it also follows from this that $\phi_*$ is
an isomorphism on $KT_2(\mathcal O_{2^k+1})$.  
Therefore these groups persist in the
direct limit and $N_2$ has the united $K$-theory as shown.  

The $K$-theory element we want $p$ to represent is the non-zero element of 
$KO_2(N_2)$, so we will consider $S^2 N_2$ to shift that $K$-theory element 
down to degree 0.  We get $N_1$ into the picture by the following 
unitization of $S^2 N_2$.  For each positive integer $n$ let
$$\phi_n \colon M_{n!}(S^2 N_2^+) 
	\rightarrow M_{(n+1)!}(S^2 N_2^+)$$  
be given by
$$\phi_n(a \oplus \lambda) = 
\begin{pmatrix}
	a \oplus \lambda & 0 & 0 & \dots & 0 \\
	0 & \lambda  & 0 & \dots & 0 \\
	0 & 0 & \lambda  & \dots & 0 \\
	\vdots & \vdots & \vdots & \ddots & \vdots \\
	0 & 0 & 0 & \dots & \lambda 
\end{pmatrix}$$
and define $N$ to be the limit of the system $\{M_{n!}(S^2 N_2^+), \phi_n \}$.

The split exact sequence
$$0 \rightarrow S^2 N_2 \rightarrow S^2 N_2^+ \rightarrow \R \rightarrow 0 \; $$
passes to a split exact sequence
$$ 0 \rightarrow \mathcal{K} \otimes S^2 N_2 
	\rightarrow N
	\rightarrow N_1
	\rightarrow 0 \; $$
in the limit, so we have $K\crt(N) \cong K\crt(N_1) \oplus K\crt(S^2 N_2)$.
According to Theorem~3.3 in \cite{Bou}, a \ct-module is injective if and only
if it is acyclic and the groups in the complex part are divisible.  Hence
$K\crt(N)$ is injective.

Let $p$ be a projection in $M_l(N)$ such that $[p] - [p_k] = 1 \oplus 1 \in 
KO_0(N) \cong KO_0(N_1) \oplus KO_2(N_2) \cong \Q \oplus \Z_2$.  Then 
$[p] = (1+k) \oplus 1$.  By picking a new preferred generator of $\Q$,
we may assume that $[p] = 1 \oplus 1$.  For future reference,
we also note that $c[p] = 1 \oplus \tfrac{1}{2} \in \Q \oplus
\Z(2^{\infty})$ and $\ve[p] = 1 \oplus 1 \in \Q \oplus \Z_2$.

Now we must show that if $F$ is a C*-algebra such that $K\crt(F)$ is free, then
the homomorphism 
$$\psi \colon K\crt(F) \rightarrow K\crt(F \otimes N)$$
induced by the inclusion map $a \mapsto a \otimes p$ is injective.  

If we define 
$\phi \colon K\crt(F) \rightarrow K\crt(F) \otimes K\crt(N)$
by 
\begin{align*}
\phi(x) &= x \otimes [p] {\text ~for~} x \in KO_*(F) \\
\phi(y) &= y \otimes c[p] {\text ~for~} y \in KU_*(F) \\
\phi(z) &= z \otimes \ve [p] {\text ~for~} y \in KT_*(F) 
\end{align*}
and we let $\alpha$ be the pairing
$K\crt(F) \otimes K\crt(N) \rightarrow K\crt(F \otimes N)\; ,$
then $\alpha \circ \phi = \psi$.
By the K\"unneth formula for united $K$-theory (Theorem~4.2 of \cite{Boer2}),
$\alpha$ is an isomorphism.  So it suffices to show that the map $\phi$ is
injective.

Since any free \ct-module can be written as a direct sum of monogenic free
\ct-modules and since the homomorphism $\phi$ respects direct sums and suspensions,
it is enough to show that $\phi$ is injective in case $F$ is one of the
three free monogenic \ct-modules $F(b,0,\R)$, $F(b,0,\C)$, and
$F(b,0,T)$, shown in Table~\ref{freemon} (or see Tables I, II, and III in 
\cite{Boer2} for the tables of operations). 

\begin{table}   
\caption{Free monogenic \ct-modules}
\label{freemon} $$\begin{array}{|c|c|c|c|c|c|c|c|c|c|}  
\hline  \hline  
n & \makebox[1cm][c]{0} & \makebox[1cm][c]{1} & 
\makebox[1cm][c]{2} & \makebox[1cm][c]{3} 
& \makebox[1cm][c]{4} & \makebox[1cm][c]{5} 
& \makebox[1cm][c]{6} & \makebox[1cm][c]{7} 
& \makebox[1cm][c]{8} \\
\hline  \hline 
F(b,0,\R)\po_n  &
\Z & \Z_2 & \Z_2 & 0 & \Z &
0 & 0 & 0 & \Z         \\
\hline  
F(b,0,\R)\pu_n   &
\Z & 0 & 
\Z & 0 & 
\Z & 0 & 
\Z & 0 & 
\Z         \\
\hline  
F(b,0,\R)\pt_n &
\Z & \Z_2 & 0 & \Z &  
\Z & \Z_2 & 0 & \Z &
\Z_k         \\
\hline \hline
F(b,0,\C)\po_n  &
\Z & 0 & \Z & 0 & \Z &
0 & \Z & 0 & \Z         \\
\hline  
F(b,0,\C)\pu_n   &
\Z^2 & 0 & 
\Z^2 & 0 & 
\Z^2 & 0 & 
\Z^2 & 0 & 
\Z^2         \\
\hline  
F(b,0,\C)\pt_n &
\Z & \Z & \Z & \Z &  
\Z & \Z & \Z & \Z &
\Z         \\
\hline \hline
F(b,0,T)\po_n &
\Z & \Z & \Z_2 & 0 & \Z & \Z &
\Z_2 & 0 & \Z          \\
\hline  
F(b,0,T)\pu_n   &
\Z & \Z & 
\Z & \Z & 
\Z & \Z & 
\Z & \Z & 
\Z         \\
\hline  
F(b,0,T)\pt_n &
\Z^2 & \Z \oplus \Z_2 & \Z_2 & \Z &  
\Z^2 & \Z \oplus \Z_2 & \Z_2 & \Z &
\Z^2         \\
\hline \hline 
\end{array}$$
\end{table}

First assume $F = F(b, 0, \R)$.  Then the homomorphism
\begin{align*}
\phi \colon K\crt(F) &\rightarrow K\crt(F) \otimes K\crt(N)   \\
	& \cong K\crt(F) \otimes K\crt(N_1)
		\oplus K\crt(F) \otimes  \Sigma^2 K\crt(N_2)   \\
	& \cong K\crt(N_1) \oplus \Sigma^2 K\crt(N_2)
\end{align*}
is the unique homomorphism which sends the generator $b$ to $b \otimes [p]
= 1 \oplus 1 \in \Q \oplus \Z_2$ and is seen to be injective by 
checking that the \ct-module generated by $1 \oplus 1$ is isomorphic 
to $F(b,0,\R)$. 

In the case of $K\crt(F) = F(b, 0, \C)$ we have the homomorphism
$$\phi \colon K\crt(F) \rightarrow 
	 K\crt(F) \otimes K\crt(N_1)
		\oplus K\crt(F) \otimes \Sigma^2 K\crt(N_2)  \; . 
$$
Using Proposition~3.6 of \cite{Boer2}, we compute $K\crt(F) \otimes
K\crt(N_1)$ and $K\crt(F) \otimes \Sigma^2 K\crt(N_2)$ as shown in
Table~\ref{F2timesN}.  The homomorphism $\phi$ sends the generator $b$ to
$b \otimes c[p] = b \otimes (1 \oplus \tfrac{1}{2})$ which is the element
$\smv{1}{0} \oplus \smv{1/2}{0} \in \Q^2 \oplus \Z(2^{\infty})^2$. 
Again, we check that the \ct-module generated by this element is 
isomorphic to $F(b,0,\C)$.

\begin{table}   
\caption{$F(b,0,\C) \otimes K\crt(N)$}  
\label{F2timesN} $$\begin{array}{|c|c|c|c|c|c|c|c|c|c|}  
\hline  \hline  
n & \makebox[1cm][c]{0} & \makebox[1cm][c]{1} & 
\makebox[1cm][c]{2} & \makebox[1cm][c]{3} 
& \makebox[1cm][c]{4} & \makebox[1cm][c]{5} 
& \makebox[1cm][c]{6} & \makebox[1cm][c]{7} 
& \makebox[1cm][c]{8} \\
\hline  \hline 
(F(b,0,\C) \otimes K\crt(N_1))\po  &
\Q & 0 & \Q & 0 & \Q &
0 & \Q & 0 & \Q         \\
\hline  
(F(b,0,\C) \otimes K\crt(N_1))\pu  &
\Q^2 & 0 & 
\Q^2 & 0 & 
\Q^2 & 0 & 
\Q^2 & 0 & 
\Q^2         \\
\hline  
(F(b,0,\C) \otimes K\crt(N_1))\pt &
\Q & \Q & \Q & \Q & \Q & \Q & \Q & \Q & \Q         \\
\hline \hline
(F(b,0,\C) \otimes K\crt(N_2))\po  &
\Z(2^{\infty}) & 0 & \Z(2^{\infty}) & 0 &  
\Z(2^{\infty}) & 0 & \Z(2^{\infty}) & 0 &
\Z(2^{\infty})         \\
\hline  
(F(b,0,\C) \otimes K\crt(N_2))\pu  &
\Z(2^{\infty})^2 & 0 & \Z(2^{\infty})^2 & 0 &  
\Z(2^{\infty})^2 & 0 & \Z(2^{\infty})^2 & 0 &
\Z(2^{\infty})^2         \\
\hline  
(F(b,0,\C) \otimes K\crt(N_2))\pt &
\Z(2^{\infty}) & \Z(2^{\infty}) & \Z(2^{\infty}) & \Z(2^{\infty}) &  
\Z(2^{\infty}) & \Z(2^{\infty}) & \Z(2^{\infty}) & \Z(2^{\infty}) &
\Z(2^{\infty})        \\
\hline \hline 
\end{array}$$
\end{table}

Finally, in the case of $K\crt(F) = F(b,0,T)$ we have the homomorphism
$$ \phi \colon K\crt(F) \rightarrow  
	K\crt(F) \otimes K\crt(N_1)
		\oplus K\crt(F) \otimes \Sigma^2 K\crt(N_2)  $$
with the latter groups shown in Table~\ref{F3timesN}, computed using
Proposition~3.7 of \cite{Boer2}.  The homomorphism $\phi$ is determined by
sending $b$ to $b \otimes \ve[p] = b \otimes (\alpha \oplus 1)$ which is
the element $\alpha \oplus 1 \in \Q^2 \oplus \Z_2$.  Again, check that 
this element generates a \ct-submodule freely.
\end{proof}

\begin{table}   
\caption{$F(b,0,T) \otimes K\crt(N)$}  
\label{F3timesN} $$\begin{array}{|c|c|c|c|c|c|c|c|c|c|}  
\hline  \hline  
n & \makebox[1cm][c]{0} & \makebox[1cm][c]{1} & 
\makebox[1cm][c]{2} & \makebox[1cm][c]{3} 
& \makebox[1cm][c]{4} & \makebox[1cm][c]{5} 
& \makebox[1cm][c]{6} & \makebox[1cm][c]{7} 
& \makebox[1cm][c]{8} \\
\hline  \hline 
(F(b,0,T) \otimes N_1)\po  &
\Q & \Q & 0 & 0 & \Q &
\Q & 0 & 0 & \Q         \\
\hline  
(F(b,0,T) \otimes N_1)\pu   &
\Q & \Q & \Q & \Q & 
\Q & \Q & \Q & \Q & 
\Q         \\
\hline  
(F(b,0,T) \otimes N_1)\pt &
\Q^2 & \Q & 0 & \Q & \Q^2 & \Q & 0 & \Q & \Q^2         \\
\hline \hline
(F(b,0,T) \otimes N_2)\po  &
\Z(2^{\infty}) & \Z(2^{\infty}) & 0 & \Z_2 &
\Z(2^{\infty}) & \Z(2^{\infty}) & 0 & \Z_2 &
\Z(2^{\infty})         \\
\hline  
(F(b,0,T) \otimes N_2)\pu  &
\Z(2^{\infty}) & \Z(2^{\infty})& \Z(2^{\infty}) & \Z(2^{\infty}) &  
\Z(2^{\infty}) & 0Z(2^{\infty}) & \Z(2^{\infty}) & \Z(2^{\infty}) &
\Z(2^{\infty})         \\
\hline  
(F(b,0,T) \otimes N_2)\pt &
\Z(2^{\infty})^2 & \Z(2^{\infty}) & \Z_2 & \Z(2^{\infty}) &  
\Z(2^{\infty})^2 & \Z(2^{\infty}) & \Z_2 & \Z(2^{\infty}) & 
\Z(2^{\infty})^2        \\
\hline \hline 
\end{array}$$
\end{table}

Let $\mathfrak{S}A = (S^+)^8 A$ be the eightfold unital suspension of
$A$ as utilized in Section~2.2 of \cite{Boer2}.

\begin{prop} \label{injectivemap}
Let $A$ be a real unital C*-algebra.  Then there exists a C*-algebra $D$ such that
$K\crt(D)$ is an injective \ct-module and a C*-algebra map 
$\nu \colon S \mathfrak{S}A \rightarrow D$ 
such that the induced map $\nu_*$ on united $K$-theory is injective.
\end{prop}

\begin{proof}
By Proposition~2.1 in \cite{Boer2}, there is exists a C*-algebra $F$ such that
$K\crt(F)$ is free and a homomorphism 
$\mu \colon F \rightarrow \mathcal{K} \otimes \mathfrak{S}A$ 
such that $\mu_*$ is surjective on united $K$-theory.  
Then as the proof of Theorem~2.2 in \cite{Boer2} (using the construction of
Schochet in section 2 of \cite{Schochet3}) let $F_0$ be the mapping cylinder of
$\mu$ and let $F_1$ be the mapping cone of $\mu$ to form the sequence
$$0 \rightarrow F_1 \xrightarrow{\mu_1} F_0 \xrightarrow{\mu_0} 
	\mathcal{K} \otimes \mathfrak{S}A \rightarrow 0  \; .$$
where $K\crt(F_i)$ is a free \ct-module for $i = 0,1$
and the induced map $(\mu_0)_*$ is surjective on united $K$-theory. 

Let $\phi \colon F_0 \rightarrow F_0 \otimes N$ be the homomorphism defined
by $\phi(a) = a \otimes p$ and form the commutative square 
$$ \xymatrix{ 
F_1 \ar[rr]^{\mu_1} \ar[d]^{=} 
&& F_0 \ar[d]^{\phi} \\
F_1 \ar[rr]^{\phi \circ \mu_1} 
&& F_0 \otimes N                        }  $$

Using the naturality of the mapping cone construction (see 
Proposition~2.9 of \cite{Schochet3}) we form a diagram 
$$ \xymatrix{
SF_1 \ar[d]^= \ar[rr]^{S\mu_1}
&&SF_0 \ar[d]^{S\phi} \ar[r]
& C\mu_1 \ar[d]^{\gamma} \ar[r]
& F_1 \ar[rr]^{\mu_1} \ar[d]^{=} 
&& F_0 \ar[d]^{\phi} \\
SF_1 \ar[rr]^{S(\phi \circ \mu_1)}
&&SF_0 \otimes N \ar[r]
& C(\phi \circ \mu_1) \ar[r]
& F_1 \ar[rr]^{\phi \circ \mu_1} 
&& F_0 \otimes N         }$$
commuting up to homotopy.  The horizontal sequences are mapping cone
sequences which induce exact sequences on united $K$-theory.  In fact, 
since $\mu_1$ and $\phi \circ \mu_1$ induce
monomorphisms on united $K$-theory, the sequences in united $K$-theory are
actually short exact sequences:   
$$ \xymatrix{
0 \ar[r] 
&K\crt(F_1) \ar[d]^= \ar[rr]^{(\mu_1)_*}
&&K\crt(F_0) \ar[d]^{\phi_*} \ar[r]
&K\crt(C\mu_1) \ar[d]^{\gamma_*} \ar[r]
&0     \\
0 \ar[r]
&K\crt(F_1) \ar[rr]^{(\phi \circ \mu_1)_*}
&&K\crt(F_0 \otimes N) \ar[r]
&K\crt(C(\phi \circ \mu_1)) \ar[r]
&0	\; .			   }$$

Now since $\phi_*$ is a monomorphism, the snake lemma (Lemma~1.3.2 in
\cite{Weibel}) shows that $\gamma_*$ is also a
monomorphism.  We will take $C(\phi \circ \mu_1)$ to be the sought-for
C*-algebra $D$.
As the quotient of an injective \ct-module, the
\ct-module $K\crt(D)$ is injective.  
Let 
$$\iota \colon S(\mathcal{K} \otimes \mathfrak{S}A) \rightarrow C\mu_1$$ 
be
the homotopy equivalence of Proposition~2.4 in \cite{Schochet3} and let
$$j \colon S \mathfrak{S}A \rightarrow S(\mathcal{K} \otimes \mathfrak{S}A)$$ 
be the homomorphism induced by the choice of a rank one
projection in $\mathcal{K}$.  We take the composition
$$\gamma \circ \iota \circ j 
\colon S \mathfrak{S}A \rightarrow C(\phi \circ \mu_1) = D$$ 
to be the sought-for homomorphism $\nu$.  Since
$\iota$ and $j$ induce isomorphisms on $K$-theory, $\nu$ induces a monomorphism
on united $K$-theory as desired. 

\end{proof}

We now commence the proof proper of Theorem~\ref{mainuct}, beginning with the 
case that $B$ is unital.

\begin{proof}[Proof of Theorem~\ref{mainuct}]

Let $A$ be a real C*-algebra whose complexification is in the complex 
bootstrap category and let $B$ be a real separable unital C*-algebra.   
Let $\nu \colon S \mathfrak{S}B \rightarrow D$ be given according to 
Proposition~\ref{injectivemap}.  Then using the natural mapping cone 
construction described in Proposition~2.5 of \cite{Schochet3} there exist 
C*-algebras $I_0$ and $I_1$ and homomorphisms $\iota_0$ and $\iota_1$ such 
that
\begin{equation} \label{georesolution}
0 \rightarrow S^2 \mathfrak{S} B \xrightarrow{\iota_0} 
	I_0 \xrightarrow{\iota_1}
	I_1 \rightarrow 0 
\end{equation}
is a semisplit (see Example~19.5.2(b) in \cite{Black}) 
exact sequence of C*-algebras.  Furthermore, $I_0$ is homotopy equivalent 
to $SD$ and under this equivalence, the 
homomorphism $\iota_0$ corresponds to $S \nu$.
Since $\nu_*$ is injective, when we apply united $K$-theory, we obtain a 
short exact sequence
\begin{equation} \label{Kresolution}
0 \rightarrow K\crt(S^2 \mathfrak{S} B) \xrightarrow{(\iota_0)_*} 
	K\crt(I_0) \xrightarrow{(\iota_1)_*} 
	K\crt(I_1) \rightarrow 0 \; .
\end{equation}
Since $K\crt(I_0)$ is an injective CRT-module and since 
$K\crt(S^2 \mathfrak{S} B)$ has injective dimension 1, we conclude that 
$K\crt(I_1)$ is injectve and we're looking at an injective resolution.

If we instead apply the functor $KK(A, \cdot)$ to Sequence~\ref{georesolution}
we obtain the sequence
$$\cdots \rightarrow 
KK\crt(A, S^2 \mathfrak{S} B) \xrightarrow{(1, \iota_0)_*} 
KK\crt(A, I_0) \xrightarrow{(1, \iota_1)_*} 
KK\crt(A, I_1) \xrightarrow{\delta} 
KK\crt(A, S^2 \mathfrak{S} B) \rightarrow \cdots$$
which we unsplice to form the short exact sequence
$$0 \rightarrow \coker(1, \iota_1)_*
	\rightarrow KK\crt(A, S^2 \mathfrak{S} B)
	\rightarrow  \ker(1, \iota_1)_*
	\rightarrow 0 \; .$$
To complete the construction of the universal coefficient exact sequence, 
it remains to make the identifications
\begin{align*}
\ker(1, \iota_1)_* &\cong \hom\scrt(K\crt(A), K\crt(SB)) \\
\coker(1, \iota_1)_* &\cong \ext\scrt(K\crt(A), K\crt(SB))
\end{align*}
which we do by use of the following commutative diagram:
\begin{equation} \label{kappadiagram}
\xymatrix{
0 \ar[d]  \\
\hom\scrt(K\crt(A), K\crt(S^2 \mathfrak{S} B))   \ar[d]^{\hom(1, (\iota_0)_*)}			
&& KK\crt(A, S^2 \mathfrak{S} B)  \ar[d]^{(1, \iota_0)_*}
		\ar[ll]_-{\gamma(A, S^2 \mathfrak{S^2}B)}   \\
\hom\scrt(K\crt(A), K\crt(I_0))  \ar[d]^{\hom(1, (\iota_1)_*)}
&& KK\crt(A, I_0)	  	 \ar[d]^{(1, \iota_1)_*} 
				 \ar[ll]_-{\gamma(A, I_0)} \\
\hom\scrt(K\crt(A), K\crt(I_1))   \ar[d]
&& KK\crt(A, I_1)	  	 \ar[d]^{\delta}
				 \ar[ll]_-{\gamma(A, I_1)}  	 \\
\ext\scrt(K\crt(A), K\crt(S^2 \mathfrak{S} B))   \ar[d]	
&& KK\crt(A, S^2 \mathfrak{S} B)		       	 \\
0	} 
\end{equation}

In this diagram, both vertical sequences are exact, 
the left one derived from the resolution given by  
Sequence~\ref{Kresolution}.  The horizontal 
homomorphisms $\gamma(A, I_0)$ and $\gamma(A, I_1)$ are isomorphisms 
since by Proposition~\ref{specialcase}.  

The diagram shows that the kernel of $(1, \iota_1)_*$ is isomorphic to 
the kernel of $\hom(1, (\iota_1)_*)$ which, in turn, is isomorphic to 
$\hom\scrt(K\crt(A), K\crt(S^2 \mathfrak{S}B))$.  
Furthermore, in this identification, the composition
$$KK\crt(A, S^2 \mathfrak{S} B) \rightarrow \ker(1, \iota)_* 
	\rightarrow \hom\scrt(K\crt(A), K\crt(S^2 \mathfrak{S} B))$$
is $\gamma(A, S^2 \mathfrak{S} B)$.

The diagram also shows that the cokernel of $(1, \iota_1)_*$ is isomorphic to 
the cokernel of $\hom(1, (\iota_1)_*)$ which is isomorphic to 
$\ext\scrt(K\crt(A), K\crt(S^2 \mathfrak{S}B))$.  Composing this 
homomorphism with $\delta$ gives the homomorphism 
$$\kappa(A, S^2 \mathfrak{S} B) \colon
\ext\scrt(K\crt(A), K\crt(S^2 \mathfrak{S}B))
\rightarrow KK\crt(A, S^2 \mathfrak{S} B)$$ 
which has degree $-1$.

This proves that the universal coefficient theorem holds for the pair $(A, 
S^2 \mathfrak{S} B)$.  In Proposition~\ref{uctnatural} we will show 
that this homomorphism
$\kappa$ is natural with respect to homomorphisms $B \rightarrow B'$.  In 
particular, consider the homomorphism $B \rightarrow 0$, which induces a 
homomorphism $\mathfrak{S} B \rightarrow \mathfrak{S} 0 = (S^+)^7 \R$.  
This homomorphism forms part of a split exact sequence
$$0 \rightarrow S^8 B \rightarrow \mathfrak{S} B \rightarrow  
\mathfrak{S} 0 \rightarrow 0 \; .$$
After suspending twice, we obtain the sequence
$$0 \rightarrow S^{10} B \xrightarrow{i} 
	S^2 \mathfrak{S} B \xrightarrow{\pi}  
	S^2 \mathfrak{S} 0 \rightarrow 0 \; .$$  
Because this sequence is split, the vertical sequences of the following 
diagram are exact:
$$\xymatrix{
& 0 \ar[d] & 0 \ar[d] & 0 \ar[d]    \\
0 \ar[r] &  
\ext\scrt(K\crt(A), K\crt(S^{10} B) )  
	 \ar[d]^{(1, i_*)}   &
KK\crt(A, S^{10} B) 
	\ar[r]^-{\gamma(A, S^{10} B)} \ar[d]^{i_*} &
\hom\scrt(K\crt(A), K\crt(S^{10} B)) 
	\ar[d]^{(1, i_*)}  \ar[r]	&
0     \\
0 \ar[r] &
\ext\scrt(K\crt(A), K\crt(S^2 \mathfrak{S} B)  
	\ar[r]^-{\kappa(A, S^2 \mathfrak{S} B)} \ar[d]^{(1, \pi_*)}   &
KK\crt(A, S^2 \mathfrak{S} B) 
	\ar[r]^-{\gamma(A, S^2 \mathfrak{S} B)} \ar[d]^{\pi_*}	&
\hom\scrt(K\crt(A), K\crt(S^2 \mathfrak{S} B)) \ar[r] \ar[d]^{(1, \pi_*)}    &
0   \\
0 \ar[r]  &
\ext\scrt(K\crt(A), K\crt(S^2 \mathfrak{S} 0)) 
	\ar[r]^-{\kappa(A, S^2 \mathfrak{S} 0)} \ar[d]	&
KK\crt(A, S^2 \mathfrak{S} 0) 
	\ar[r]^-{\gamma(A, S^2 \mathfrak{S} 0)} \ar[d]	&
\hom\scrt(K\crt(A), K\crt(S^2 \mathfrak{S} 0))
	\ar[r] \ar[d] &
0    \\
& 0 & 0 & 0
} $$
Furthermore, we have shown that the two lower horizontal sequences are 
exact, the map $\gamma$ commutes with $i_*$ and with $\pi_*$, and the 
map $\kappa$ commutes with $\pi_*$ by the proof of 
Proposition~\ref{uctnatural}.  Therefore, a homomorphism 
$$\kappa(A, S^{10}) \colon \ext\scrt(K\crt(A), K\crt(S^{10} B))
	\rightarrow KK\crt(A, S^{10} B)$$
can be defined to make the diagram commute and then by a diagram chase
(Exercise~1.3.2 of \cite{Weibel}), the top horizontal sequence is exact.  
This proves that the 
universal coefficient theorem holds for the pair $(A, S^{10} B)$.  By 
Lemma~\ref{suspension} below, it holds for the pair $(A, B)$.

Finally, an argument similar to that in the previous paragraph, based on 
the exact sequence
$$0 \rightarrow B \rightarrow B^+ \rightarrow \R \rightarrow 0$$
takes care of the general case in which $B$ is not necessarily unital.
\end{proof}

\begin{lemma} \label{suspension}
The universal coefficient theorem holds for a pair $(A,B)$ if and only if 
it holds for the pair $(A, SB)$.
\end{lemma}

\begin{proof}
Consider the following diagram.
$$\xymatrix{
0 \ar[d] & & 0 \ar[d]   \\
\ext\scrt(K\crt(A), K\crt(B)) \ar[d]^{\kappa(A, B)} \ar[rr]     & &
\ext\scrt(K\crt(A), K\crt(SB)) \ar[d]^{\kappa(A, SB)}  \\
KK\crt(A, B)    \ar[d]^{\gamma(A, B)} \ar[rr]    & &
KK\crt(A, SB)   \ar[d]^{\gamma(A, SB)}     \\
\hom\scrt(K\crt(A), K\crt(B))   \ar[d]  \ar[rr]   & &
\hom\scrt(K\crt(A), K\crt(SB))   \ar[d]   \\
0 & &  0  } $$
The horizontal homomorphisms are isomorphisms of degree $-1$.  Therefore, 
there is a homomorphism $\kappa(A, B)$ as shown in the diagram making the 
left vertical sequence exact if and only if there is a homomorphism 
$\kappa(A, SB)$ making the right vertical sequence exact.
\end{proof}

\begin{prop} \label{uctnatural}
The universal coefficient sequence is natural with respect to homomorphisms of $A$ and $B$.
\end{prop}

\begin{proof}
It is clear that $\gamma$ is natural.  We must show that $\kappa$ is 
natural in both arguments.  First we show that $\kappa$ doesn't depend on 
the particular resolution of $S^2 \mathfrak{S} B$ chosen in its 
construction.  For this suppose that we start with two alternatives
$$\nu_1 \colon S \mathfrak{S} B \rightarrow D_1$$
and
$$\nu_2 \colon S \mathfrak{S} B \rightarrow D_2$$
which satisfy the requirements described in Proposition~\ref{injectivemap}.  Then we 
construct yet another alternative
$$\nu_3 \colon S \mathfrak{S} B \rightarrow D_3$$
where $D_3 = D_1 \oplus D_2$ and $\nu_3 = \nu_1 \oplus \nu_2$.  
These three possibilities lead to three possible injective resolutions
$$0 \rightarrow S^2 \mathfrak{S} B \xrightarrow{\iota_{i,0}} 
		I_{i,0} \xrightarrow{\iota_{i,1}} 
		I_{i,1} \rightarrow 0$$
for $i \in \{1,2,3\}$ which lead to three possible constructions of 
$\kappa$.  Furthermore, for each $i \in \{1,2\}$, the coordinate projection $D_3 
\rightarrow D_i$ commutes with the $\nu$ homomorphisms, and thus induces 
a homomorphism of injective resolutions:
$$\xymatrix{
0 \ar[r]
& S^2 \mathfrak{S} B \ar[r]^{\iota_{3,0}} \ar[d]^{=}
& I_{3,0} \ar[r]^{\iota{3,1}}                \ar[d] 
& I_{3,1} \ar[r] \ar[d]
& 0  \\
0 \ar[r]
& S^2 \mathfrak{S} B \ar[r]^{\iota_{i,0}}
& I_{i,0} \ar[r]^{\iota_{i,1}}
& I_{i,1} \ar[r]
& 0 			} $$

In turn, this induces a homomorphism from the groups of Diagram~\ref{kappadiagram} 
based on the resolution built from $D_3$ to the corresponding groups of the 
same diagram based on the resolution built from $D_i$.  The resulting three 
dimensional diagram 
(which we do not attempt to represent) commutes and shows that the 
homomorphism $\kappa$ is the same in either case. 

Now if $\phi \colon A \rightarrow A'$ is a homomorphism of real 
C*-algeberas, then we choose a fixed geometric injective resolution of 
$S^2 \mathfrak{S} B$ to compute $\kappa(A, S^2 \mathfrak{S} B)$ and 
$\kappa(A', S^2 \mathfrak{S} B)$.  Because $\gamma$ is natural, there is a 
homomorphism induced by $\phi$ from the groups of 
Diagram~\ref{kappadiagram} to the groups of the same diagram with $A$ 
replaced by $A'$.  Since $\kappa$ is defined via this diagram in terms of 
$\gamma$, it follows that $\kappa$ commutes with $\phi$.

Finally suppose that $\phi \colon B \rightarrow B'$ is a homomorphism of 
real C*-algebras.  First construct two homomorphisms
$$\nu \colon S \mathfrak{S} B \rightarrow D$$
and
$$\nu' \colon S \mathfrak{S} B' \rightarrow D'$$
according to Proposition~\ref{injectivemap}.  Then consider the commutative 
diagram
$$\xymatrix{
S \mathfrak{S} B   \ar[r]^{\nu''}  \ar[d]^{\phi}
& D \oplus D'   \ar[d]^{\smh{0}{1}} \\
S \mathfrak{S} B' \ar[r]^{\nu'}
& D'
}$$
where $\nu'' = \nu \oplus \nu' \circ S \mathfrak{S} B \phi$.  
Again this gives us a homomorphism from one geometric injective 
resolution to another and thus we obtain a family of homomorphisms from 
the version of Diagram~\ref{kappadiagram} based on the modified resolution 
of $S^2 \mathfrak{S} B$ to another the version of 
Diagram~\ref{kappadiagram} based on the resolution of $S^2 \mathfrak{S} 
B'$.  Since $\kappa$ does not depend on the particular geometric resolution, 
it follows that $\kappa$ commutes with $\phi$.
\end{proof}

\subsection{Some Corollaries} \label{corsection}

If we focus on the real part of each \ct-module
in the short exact sequence of Theorem~\ref{mainuct},
we obtain the following corollary.  Recall 
that $[M,N]$ represents the graded group of \ct-morphisms and let 
$\ext_{[~,~]}(M,N)$ denote the associated derived functor.  

\begin{cor} \label{realuct}
Let $A$ and $B$ be real separable C*-algebras 
such that $\C \otimes A$ is in $\mathcal{N}$.  
Then there is a short
exact sequence
\end{cor}
$$0 \rightarrow {\rm Ext}_{[~,~]}(K\crt(A), K\crt(B)) \xrightarrow{\kappa\so}
	KK(A,B) \xrightarrow{\gamma\so}
	[K\crt(A), K\crt(B)] 
	\rightarrow 0$$
where $\kappa\so$ has degree $-1$ and the homomorphism $\gamma\so$ is given 
by $\gamma\so(x)(y) = \alpha\so(y \otimes x)$.

\begin{proof}
By definition, $KK\crt(A,B)\po = KK(A,B)$ and by 
Proposition~\ref{homo=[]}, 
$$\hom\scrt(K\crt(A), K\crt(B))\po \cong [K\crt(A), K\crt(B)] \; .$$
To identify the real part of $\gamma$ under these 
identifications, let $x \in KK_i(A,B)$.  Then 
$\gamma(x) \in \hom\scrt(K\crt(A), K\crt(B))\po
	= [K\crt(A) \otimes\scrt F(b,i,\R), K\crt(B)]$
is defined by $\gamma(x) = \Gamma(\alpha)(x) = \alpha \circ (1 \otimes 
\mu_b^x)$ (see the proof to Proposition~\ref{adjoint} for the description 
of $\Gamma$) thus $\gamma\so(x)(y) = \alpha\so(1 \otimes \mu_b^x)(y \otimes b) = 
\alpha\so(y \otimes x)$.

Finally, we claim that $\ext\scrt(K\crt(A), K\crt(B))\po 
\cong \ext_{[~,~]}(K\crt(A), K\crt(B))$.  Let 
$$0 \rightarrow F_1 \xrightarrow{\mu_1} F_0 \xrightarrow{\mu_0} K\crt(A) 
\rightarrow 0$$
be a free resolution of $K\crt(A)$.  We apply $\hom\scrt$ to obtain
\begin{multline*}
0 \leftarrow \ext\scrt(K\crt(A), K\crt(B))
\leftarrow \hom\scrt(F_1, K\crt(B)) \xleftarrow{\mu_1^*}   \\
\hom\scrt(F_0, K\crt(B)) \xleftarrow{\mu_0^*}
\hom\scrt(K\crt(A), K\crt(B)) \leftarrow 0
\end{multline*}
and then restrict to the real part to obtain
$$0 \leftarrow \ext\scrt(K\crt(A), K\crt(B))\po
\leftarrow [F_1, K\crt(B)] \xleftarrow{\mu_1^*}
[F_0, K\crt(B)] \xleftarrow{\mu_0^*}
[K\crt(A), K\crt(B)] \leftarrow 0 \; .$$
But the cokernel of $\mu_1^*$ is by definition 
$\ext_{[~,~]}(K\crt(A), K\crt(B))$, proving the claim.
\end{proof}

The final corollaries can be proven from Corollary~\ref{realuct} exactly as 
in the proof of Propositions~23.10.1 and 23.10.2 in \cite{Black}.

\begin{cor}
Let $A$ and $B$ be real separable C*-algebras 
such that $\C \otimes A$ is in $\mathcal{N}$, and let $x \in KK(A,B)$ have 
the property that $\gamma(x) \in [K\crt(A), K\crt(B)]$ is an isomorphism.  
Then $x$ is a $KK$-equivalence.
\end{cor} 

\begin{cor}
Let $A$ and $B$ be real separable C*-algebras such that $\C \otimes A$ 
and $\C \otimes B$ are in $\mathcal{N}$.  Then $A$ and $B$ are 
$KK$-equivalent if and only if $K\crt(A) \cong K\crt(B)$.
\end{cor}



\newpage

\vspace*{-0.2in}

\begin{thebibliography}{99}
\raggedbottom
\itemsep12pt

\bibitem{Ati2}
{\scshape M.F. Atiyah}, Vector Bundles and the K\"unneth Formula,  {\it Topology} 1 
(1962), 245--248. 

\bibitem{Black}
{\scshape B. Blackadar}, {\it $K$-theory for Operator Algebras, Second Edition},
Mathematical Sciences Research Institute Publications 5, Cambridge University Press,
1998.

\bibitem{Boer}
{\scshape J. L. Boersema}, Real C*-Algebras, United $K$-Theory, and the K\"unneth 
Formula, Ph.D. Dissertation, University of Oregon, Eugene (1999).

\bibitem{Boer2}
{\scshape J. L. Boersema}, Real C*-Algebras, United $K$-Theory, and the K\"unneth
Formula, {\it $K$-theory} 26 (2002), 345--402.

\bibitem{Bou} 
{\scshape A.K. Bousfield}, A Classification of $K$-Local Spectra,  {\it Journal of 
Pure and Applied Algebra } 66 (1990), 121--163.

\bibitem{Elliot}
{\scshape G. Elliott}, On the classification of inductive limits of 
sequences of semisimple finite-dimensional algebras, {\it J. Algebra} 
38 (1976), 29--44. 

\bibitem{EG}
{\scshape G. Elliott and G. Gong}, On the classification of C*-algebras of 
real rank zero, II, {\it Annals of Mathematics}
144 (1996), 497--610.

\bibitem{Gio}
{\scshape T. Giordano}, A classification of approximately finite real 
C*-algebras, {\it J. reine angew. Math.} 
385 (1988), 161--194. 

\bibitem{Goodearl}
{\scshape K.R. Goodearl}, {\it Notes on Real and Complex C*-algebras} Shiva Publishing 
Limited, 1982. 

\bibitem{Hewitt}
{\scshape B. Hewitt}, On the homotopical classification of $KO$-module 
spectra,  Ph.D. Dissertation,  University 
of Illinois, Chicago (1996).

\bibitem{JT}
{\scshape K. Jensen and K. Thomsen}, {\it Elements of $KK$-Theory}
Birkh\"auser, 1991.

\bibitem{Kasp}
{\scshape G.G. Kasparov}, The Operator $K$-Functor and Extensions of C*-Algebras, {\it 
Math. USSR Izvestija} 16 (1981), 513--572.

\bibitem{Kirchberg}
{\scshape E. Kirchberg}, The classification of purely infinite C*-algebras using 
Kasparov's theory, in preparation.

\bibitem{Lin}
{\scshape H. Lin}, Classification of simple C*-algebras of tracial 
topological rank zero, preprint.

\bibitem{PS}
{\scshape W. Paschke and N. Salinas}, Matrix algebras over $\mathcal{O}_n$, 
{\it Michigan Journal of Mathematics} 26 (1979), 3--12.

\bibitem{Palmer}
{\scshape T.W. Palmer}, Real C*-Algebras, {\it Pacific Journal of Mathematics} 35, 
No. 1, 1970.

\bibitem{Phillips}
{\scshape N.C. Phillips}, A classification theorem for nuclear purely infinite simple 
C*-algebras, {\it Documenta Mathematica} 5 (2000), 49--114.

\bibitem{RS}
{\scshape J. Rosenberg and C. Schochet}, The K\"unneth thoerem and the universal 
coefficient theorem for Kasparov's generalized $K$-functor, {\it Duke Mathematical
Journal} 55 (1987), 431--474.

\bibitem{Rordam}
{ M. R\o rdam}, ``Classification of Nuclear C*-Algebras,'' 
Springer-Verlag, 2002.

\bibitem{Schochet}
{\scshape C. Schochet}, Topological Methods for C*-Algebras II:  Geometric Resolutions 
and the K\"unneth Formula, {\it Pacific Journal of Mathematics} 98 (1982), 
443--458.

\bibitem{Schochet3}
{\scshape C. Schochet}, Topological Methods for C*-Algebras III:  Axiomatic Homology, 
{\it Pacific Journal of Mathematics} 114 (1984), 399--445.

\bibitem{Schroder}
{\scshape H. Schr\"oder}, {\it $K$-Theory for real C*-algebras and applications} 
Pitman Research Notes in Mathematics Series 290, 1993.

\bibitem{Weibel}
{\scshape C.A. Weibel}, ``An Introduction to Homological Algebra,'' Cambridge
studies in advanced mathematics 38, Cambridge University Press, 1994.

\end{thebibliography}
\end{document}